\documentclass[amstex,12pt, amssymb]{article}

\usepackage{mathtext}
\usepackage[cp1251]{inputenc}
\usepackage[T2A]{fontenc}
\usepackage[dvips]{graphicx}
\usepackage{amsmath}
\usepackage{amssymb}
\usepackage{amsxtra}
\usepackage{latexsym}
\usepackage{ifthen}

\newcounter{lemma}[section]

\newcounter{corollary}[section]

\newcounter{remark}[section]

\newcounter{theorem}[section]

\newcounter{proposition}[section]

\newcounter{example}

\numberwithin{equation}{section}

\begin{document}

\markboth{\centerline{E.~SEVOST'YANOV}}{\centerline{ON THE BOUNDARY
AND GLOBAL BEHAVIOR...}}

\def\cc{\setcounter{equation}{0}
\setcounter{figure}{0}\setcounter{table}{0}}

\overfullrule=0pt

%\normalsize\large

\author{EVGENY SEVOST'YANOV}

\title{
{\bf ON THE BOUNDARY AND GLOBAL BEHAVIOR OF MAPPINGS OF RIEMANNIAN
SURFACES}}

\date{\today}
\maketitle

%\large
\begin{abstract}
In this article, we study non-homeomorphic mappings of Riemannian
surfaces of the Sobolev class. We have established estimates for the
distortion of the modulus of families of paths, and as a
consequence, we obtained results on the boundary behavior of such
mappings between domains of Riemannian surfaces.
\end{abstract}

\bigskip
{\bf 2010 Mathematics Subject Classification: Primary 30C65;
Secondary 32U20, 31B15}

\section{Introduction}

One of the main problems of modern analysis is the extension of
mappings to the boundary of a domain. There are a number of results
on this topic related, in particular, to the theory of
quasiconformal mappings and their generalizations,
see~\cite{Na$_1$}--\cite{RV}. Among the listed papers, we note the
fundamental assertion of N\"{a}kki, see~\cite[Theorem~2.4]{Na$_1$}
(see also~\cite[Theorem~17.15]{Va})

\medskip
{\bf Theorem (on the extension of quasiconformal mappings to the
boundary)}. {\sl Let $D, D^{\,\prime}$ be domains in ${\Bbb R}^n,$
$n\geqslant 2,$ and let $f$ be a quasiconformal mapping of $D$ onto
$D^{\,\prime}.$ Suppose that $D$ is locally connected on its
boundary, and $\partial D^{\,\prime}$ is quasiconformally
accessible. Then $f$ has a continuous extension
$\overline{f}:\overline{D}\rightarrow\overline{{\Bbb R}^n}.$}

\medskip
This result was developed in a number of papers by other authors.
First of all, Srebro and Vuorinen extended N\"{a}kki's theorem to
quasiregular mappings, see~\cite[Theorem~4.2]{Sr} and
\cite[Theorem~4.10.II]{Vu$_1$}. More recently, Martio, Ryazanov,
Srebro and Yakubov, as well as Ignat'ev and Ryazanov obtained
results on the boundary extension of homeomorphisms with unbounded
characteristic, see~\cite[Lemma~5.16]{MRSY$_1$} and \cite[Lemma~2.1,
Corollary~2.1]{IR$_1$}. Later they were carried over to Riemannian
manifolds and metric measure spaces (see, e.g.~\cite{RSa},
\cite{Sm}, \cite{ARS} and \cite{Sev}).

\medskip
Let us now dwell on the recent results of Ryazanov and
Volkov~\cite{RV}. Here the authors proposed an approach that allows
one to study the boundary behavior of Sobolev classes acting between
two Riemannian surfaces. It should be noted that the paper~\cite{RV}
contains important results in this direction, however, they all
concern only homeomorphisms. In our opinion, it would be important
to describe the boundary behavior of similar mappings with
branching, and this problem is largely solved in this article. As
in~\cite{RV}, the main research tool is the modulus method. To a
large extent, our publication is devoted to the development of the
modulus method and the identification of fundamental opportunities
of this method in this context.

Conventionally, the paper may be divided into three parts:
establishing an estimate for the distortion of the modulus of
families of curves, $\S\S $~\ref{sec2} and \ref{sec3}; boundary
behavior of mappings, $\S~\ref{sec4}$; local the behavior of the
mappings (at the inner points of the domain) and the global behavior
of the mappings (at the inner and boundary points of the domain),
$\S \S$~\ref{sec5}--\ref{sec6}. In the last sections, homeomorphisms
and mappings with branching are studied separately, since the
corresponding results have been proved under various conditions and
are new in both cases.

Here are some definitions. A {\it Riemannian surface} is a
two-dimensional manifold with a countable base in which transition
maps between corresponding maps are conformal, see, e.g.,~\cite{RV}.
The Riemannian surface ${\Bbb S}$ considered below will be assumed
to be a {\it of hyperbolic type}, that is, a surface conformally
equivalent to the unit disk ${\Bbb D}=\{z\in {\Bbb C}: |z|<1\}$ with
''glued'' points (see \cite[\S\,6, Section~1]{KAG}). In other words,
we consider those and only those Riemannian surfaces that are
conformally equivalent to the factor space ${\Bbb D}/G, $ where $G$
is some group of linear fractional automorphisms of the unit disk
that has no fixed points and acts discontinuously in ${\Bbb D}.$
Recall that each element $p_0$ of the factor space ${\Bbb D}/G$ is
an {\it orbit} of the point $z_0\in {\Bbb D},$ that is, $p_0=\{z\in
{\Bbb D}: z=g(z_0), g\in G\}. $ In what follows we identify the
Riemannian surface $ {\Bbb S}$ with its factor representation ${\Bbb
D}/G$ in terms of some group $G$ of linear fractional mappings
$g:{\Bbb D}\rightarrow {\Bbb D}.$

In what follows, we use the {\it hyperbolic metric} on the unit disk
${\Bbb D}$ defined by the equalities
\begin{equation}\label{eq3}
h(z_1, z_2)=\log\,\frac{1+t}{1-t}\,,\quad
t=\frac{|z_1-z_2|}{|1-z_1\overline{z_2}|}\,,
\end{equation}
as well as the {\it hyperbolic area} of the set $S\subset {\Bbb D}$
and the {\it length} of the path $\gamma: [a, b]\rightarrow {\Bbb
D},$ which are given, respectively, by the relations
\begin{equation}\label{eq1}
h(S)=\int\limits_S\frac{4\,dm(z)}{(1-|z|^2)^2}\,,\quad
s_h(\gamma)=\int\limits_{\gamma}\frac{2\,|dz|}{1-|z|^2}\,, \quad
z=x+iy
\end{equation}
(see \cite[(2.4), (2.5)]{RV}). It is easy to verify by direct
calculations that the hyperbolic metric, length, and area are
invariant under linear fractional mappings of the unit disk onto
itself.

\medskip
In what follows, for a point $y_0\in {\Bbb D}$ and a number
$r\geqslant 0,$ we define the {\it hyperbolic disk} $B_h(y_0, r)$
and the {\it hyperbolic circle} $S_h(y_0, r)$ by means of equalities
\begin{equation}\label{eq7}
B_h(y_0, r):=\{y\in {\Bbb D}: h(y_0, y)<r\}\,, S_h(y_0, r):=\{y\in
{\Bbb D}: h(y_0, y)=r\}\,. \end{equation}
We perform the metrization of the surface ${\Bbb D}/G$ as follows.
If $p_1, p_2 \in {\Bbb D}/ G,$ we put
\begin{equation}\label{eq2}
\widetilde{h}(p_1, p_2):=\inf\limits_{g_1, g_2\in G}h(g_1(z_1),
g_2(z_2))\,,
\end{equation}
where $p_i=G_{z_i}=\{\xi\in {\Bbb D}:\,\exists\, g\in G:
\xi=g(z_i)\},$ $i=1,2.$ In the latter case, the set $G_{z_i}$ will
be called the {\it orbit} of the point $z_i,$ and $p_1$ and $p_2$
will be called the {\it orbits} of the points $z_1$ and $z_2,$
respectively. Note that $\widetilde{h}$ is a metric on ${\Bbb D}/G$
(see Section~2 in~\cite{RV}).

Everywhere below, ${\Bbb S}={\Bbb D}/G$ and ${\Bbb S}_*={\Bbb
D}/G_*$ are two different Riemannian surfaces of a hyperbolic type.
In what follows, we do not distinguish between the original
Riemannian surface ${\Bbb S}$ and its factor representation ${\Bbb
D}/G.$ A continuous extension of the mapping $f:D\rightarrow D_*$ to
a point $f:D\rightarrow D_*$ $p_0\in \overline{D},$ as well as other
concepts related to limit, continuity, topology, etc., are
understood in the sense of the metric spaces $({\Bbb D}/G,
\widetilde{h})$ and $({\Bbb D}/G_*, \widetilde{h}_*),$ where
$\widetilde{h}$ and $\widetilde{h}_*$ are metrics defined
in~(\ref{eq2}). The elements of length and area on the surfaces
${\Bbb S}$ and ${\Bbb S}_*$ are denoted $ds_{\widetilde{h}},$
$d\widetilde{h}$ and $ds_{\widetilde{h_*}},$ $d\widetilde{h_*},$
respectively.

\medskip
Let $D$ and $D_{\,*}$ be domains on Riemannian surfaces ${\Bbb S}$
and ${\Bbb S}_{\,*},$ respectively. A mapping $f:D\rightarrow
D_{\,*}$ is called {\it discrete} if the preimage $f^{-1}(y)$ of any
point $y \in D_{\,*}$ consists of isolated points only. A mapping
$f:D\rightarrow D_{\,*}$ is called {\it open} if the image of any
open set $U \subset D $ is an open set in $D_{\,*}.$ The definition
of mappings of the Sobolev class $W_{\rm loc}^{1,1}$ on a Riemannian
surface can be found, for example, in~\cite{RV}. In further, for
mappings $f:D\rightarrow D_{\,*}$ of a class $W_{\rm loc}^{1,1}$ in
local coordinates $f_{\overline{z}} = \left(f_x + if_y\right)/2$ and
$f_z = \left(f_x - if_y\right)/2,$ $z=x+iy.$  In addition, the {\it
norm} and the {\it Jacobian} of the mapping of $f$ in local
coordinates are expressed, respectively, by the equalities $\Vert
f^{\,\prime}(z)\Vert=|f_z|+|f_{\overline{z}}|$ and
$J_f(z)=|f_z|^2-|f_{\overline{z}}|^2.$ A {\it dilatation} of the
mapping $f$ at a point $z$ is defined by the relation
\begin{equation}\label{eq16}
K_f(z)=\frac{|f_z|+|f_{\overline{z}}|}{|f_z|-|f_{\overline{z}}|}
\end{equation}
for $J_f(z)\ne 0,$ $K_f(z)=1$ for $\Vert f^{\,\prime}(z)\Vert=0$ and
$K_f(z)=\infty$ otherwise. It is easy to verify by direct
calculations that $K_f (z)$ does not depend on local coordinates. A
mapping $f:D\rightarrow D_{\,*}$ is called a {\it mapping with
finite distortion,} if $f\in W_{\rm loc}^{1, 1}(D)$ and, in
addition, there is almost everywhere a finite function $K(z)$ such
that $\Vert f^{\,\prime}(z)\Vert^2\leqslant K(z)\cdot J_f(z)$ for
almost all $z \in D.$

\medskip
As usual, a path $\gamma$ on the Riemannian surface ${\Bbb S}$ is
defined as a continuous mapping $\gamma:I\rightarrow {\Bbb S},$
where $I$ is a finite segment, an interval or a half-interval of a
real axis. Let $\Gamma$ be a family of paths in ${\Bbb S}.$ A Borel
function $\rho:{\Bbb S}\rightarrow [0, \infty]$ is called {\it
admissible} for the family $\Gamma$ of paths $\gamma,$ if
$\int\limits_{\gamma}\rho(p)\,ds_{\widetilde{h}}(p)\geqslant 1$ for
any path $\gamma \in \Gamma.$ The latter is briefly written in the
form: $\rho\in {\rm adm}\,\Gamma.$  A {\it modulus} of the family
$\Gamma$ is a real-valued function
$$M(\Gamma):=\inf\limits_{\rho\in {\rm adm}\,\Gamma}\int\limits_{\Bbb
S}\rho^2(p)\,d\widetilde{h}(p)\,.$$
To avoid confusion, we introduce separately the notation for the
modulus of the family $\Gamma$ in the Euclidean sense, namely, put
$$M_e(\Gamma):=\inf\limits_{\rho\in {\rm
adm}_e\,\Gamma}\int\limits_{\Bbb D}\rho^2(z)\,dm(z)\,,$$
where $\rho\in {\rm adm}_e\,\Gamma$ if and only if
$\int\limits_{\gamma}\rho(z)\,|dz|\geqslant 1$ for any (locally
rectifiable) path $\gamma \in \Gamma.$

\medskip
Let $D$ and $D_{\,*}$ be given domains lying in the Riemannian
surfaces ${\Bbb S}$ and ${\Bbb S}_*,$ respectively. Given a mapping
$f:D\,\rightarrow\,D_{\,*},$ a set $E\subset D$ and
$y\,\in\,D_{\,*},$ we define a {\it multiplicity function} $N(y, f,
E)$ as the number of preimages of the point $y$ in $E,$ that is
\begin{equation}\label{eq17}
N(y,f,E)\,=\,{\rm card}\,\left\{p\in E: f(p)=y\right\}\,, \quad
N(f,E)\,=\,\sup\limits_{y\in{\Bbb S}_*}\,N(y,f,E)\,.
\end{equation}
We say that the function $\rho\colon{\Bbb S}\rightarrow[0, \infty]$
measurable with respect to the area $\widetilde{h}$ is {\it
extensively admissible} for the family $\Gamma,$ abbr.~$\rho\in{\rm
ext}\,{\rm adm}\,\Gamma,$ if the inequality
$\int\limits_{\gamma}\rho\,ds_{\widetilde{h}}(p)\geqslant 1$ is
satisfied for all locally rectifiable paths $\gamma\in
\Gamma\setminus\Gamma_0,$ where $M(\Gamma_0)=0.$

\medskip
The next class of mappings is related to the Gehring ring definition
of quasiconformality (see~\cite{Ge}) and is necessary from a
technical point of view. Its definition includes a distortion of the
modulus of families of paths in such a way as is necessary to solve
the corresponding problems of interest to us. Let $D$ and $D_{\,*}$
be domains lying in the Riemannian surfaces ${\Bbb S}$ and ${\Bbb
S}_*,$ respectively, and $Q\colon D\rightarrow(0,\infty)$ be a
measurable function with respect to the measure  $\widetilde{h}$ on
${\Bbb S}.$ We say that $f\colon D\rightarrow D_{\,*}$ is a {\it
lower $Q$-mapping} at a point $p_0\in \overline{D},$ if for some
$\varepsilon_0=\varepsilon_0(p_0)>0,$
$\varepsilon_0<d_0=\sup\limits_{p\in D}\widetilde{h}(p, p_0),$ and
any ring $\widetilde{A}(p_0, \varepsilon, \varepsilon_0)=\{p\in
{\Bbb S}: \varepsilon<\widetilde{h}(p, p_0)<\varepsilon_0\}$ the
inequality
\begin{equation}\label{eq1A}
M(f(\Sigma_{\varepsilon}))\geqslant \inf\limits_{\rho\in{\rm
ext\,adm}\,\Sigma_{\varepsilon}}\int\limits_{D\cap
\widetilde{A}(p_0, \varepsilon,
\varepsilon_0)}\frac{\rho^2(p)}{Q(p)}\,d\widetilde{h}(p)
\end{equation}
holds, where $\Sigma_{\varepsilon}$ denotes the family of all
intersections of circles $\widetilde{S}(p_0, r)=\{p\in {\Bbb S}:
\widetilde{h}(p, p_0)=r\}$ with $D,$ $r\in (\varepsilon,
\varepsilon_0)$ (see~\cite[Chapter~9]{MRSY}).

\medskip
The next assertion contains a fundamental estimate for the
distortion of families of paths in Sobolev classes (see
also~\cite[Lemma~3.1]{RV} and~\cite[Lemma~3.1]{LSS} in this regard).
Further studies related to the boundary behavior and equicontinuity
of mappings are based on estimates of such a plan (see
paragraphs~\ref{sec4}--\ref{sec6}).

\medskip
\begin{theorem}\label{thOS4.2}{\sl\, Let $D$ and $D_{\,*}$ be
domains of Riemannian surfaces ${\Bbb S}$ and ${\Bbb S}_*,$
respectively, $\overline{D_{\,*}}$ is compact in ${\Bbb S}_*$ and
$p_0\in\overline{D}.$ Then any open discrete mapping $f:D\rightarrow
D_{\,*}$ of finite distortion of the class $W^{1, 1}_{\rm loc}$ such
that $N(f, D)<\infty$ satisfies~(\ref{eq1A}) at $p_0$ for
$Q(p)=c\cdot N(f, D)\cdot K_f(p),$ where $K_f(p)$ is defined by the
relation~(\ref{eq16}), the function $N(f, D)$ is given in
(\ref{eq17}), and $c>0$ is some constant depending only on the point
$p_0$ and the domain $D_{\,*}.$ }
\end{theorem}

\medskip
Let us now formulate the main result of the article concerning the
continuous boundary extension of mappings. Let $E,$ $F\subset {\Bbb
S}$ be arbitrary sets. In the future, everywhere by $\Gamma(E, F,
D)$ we denote the family of all paths $\gamma:[a,b]\rightarrow D,$
which join $E$ and $F$ in $D,$ that is,  $\gamma(a)\in E,$
$\gamma(b) \in F$ and $\gamma(t)\in D$ for $t\in (a, \, b). $ Let us
agree to say that the boundary $\partial G$ of the domain $G$ is
{\it strongly accessible at the point $p_0\in \partial G $} if for
each neighborhood $U$ of $p_0$ there is a compactum $E\subset G,$ a
neighborhood $V\subset U $ of the same point and a number $\delta>0$
such that $M(\Gamma(E, F, G)) \geqslant \delta$ for any continua $E$
and $F$ intersecting both $\partial U,$ and $\partial V.$ We will
also say that a boundary $\partial G$ is {\it strongly accessible }
if it is strongly accessible at each of its points. Note that, for
an open closed mapping $f: D\rightarrow D_{\,*},$ the following
condition holds: $N (f, D)<\infty$ (see~\cite[Theorem~5.5]{Va$_1$}).

\medskip
A mapping $f:D\rightarrow {\Bbb C},$ $D\subset {\Bbb C}, $ is called
{\it quasiconformal} if $f$ is a homeomorphism and, moreover, there
is a constant $K\geqslant 1 $ such that $M_e(f(\Gamma))\leqslant K
\cdot M_e(\Gamma)$ for any family of paths $\Gamma$ in $D.$ We say
that the boundary of a domain $D$ in ${\Bbb S}$ is {\it locally
quasiconformal} if each point $p_0\in \partial D $ has a
neighborhood $U$ in ${\Bbb S},$ which can be mapped by a
quasiconformal mapping $\psi$ onto the unit disk ${\Bbb D}\subset
{\Bbb C}$ so that $\psi(\partial D \cap U)$ is the intersection of
${\Bbb D}$ with the straight line $x=0,$ where $z= x + iy \in {\Bbb
D}.$ The most important result of this article can be formulated as
follows.

\begin{theorem}\label{th2}{\sl\, Let $D$ and $D_{\,*}$
be domains on the Riemannian surfaces ${\Bbb S}$ and ${\Bbb S}_*,$
respectively, let $\overline{D_{\,* }}$ be a compactum in ${\Bbb
S}_*,$ let $p_0\in \partial D $ and let $Q:{\Bbb S} \rightarrow (0,
\infty)$ be a measurable function with respect to the measure
$\widetilde{h},$ $Q(p)\equiv 0$ in ${\Bbb S} \setminus D.$ Let also
$f:D\rightarrow D_{\,*}$ be an open discrete closed $W^{1, 1}_{\rm
loc}$-mapping with a finite distortion of $D$ onto $D_*$ such that
$K_f (p)\leqslant Q(p)$ for almost all $p\in D.$ Suppose that the
domain $D$ has a locally quasiconformal boundary, and the boundary
of the domain $D_{\, *}$ is strongly accessible. If the relations
\begin{equation}\label{eq45}
\int\limits_{\varepsilon}^{\varepsilon_0} \frac{dt}{\Vert
Q\Vert(t)}<\infty\,,\qquad \int\limits_0^{\varepsilon_0}
\frac{dt}{\Vert Q\Vert(t)}=\infty\,,
\end{equation}
hold for some $\varepsilon_0>0$ and any $\varepsilon\in (0,
\varepsilon_0),$ then $f$ has a continuous extension to~$p_0.$ Here
$$\Vert Q\Vert(r):=\int\limits_{\widetilde{S}(p_0,
r)}Q(p)\,\,ds_{\widetilde{h}}(p)$$ denotes $L_1$-norm of the
function $Q$ over the circle $\widetilde{S}(p_0, r).$}
\end{theorem}

\section{Preliminaries}\label{sec2}

Let us start the section with the next most important Remark.

\begin{remark}\label{rem1} Following~\cite[\S\,7.2]{Berd},
the hyperbolic distance $h(z_1, z_2) $ in the unit
disk~\cite[\S\,7.2]{Berd}, ${\Bbb D}$ can be equivalently is defined
as $\inf $ of hyperbolic lengths of all piecewise smooth paths
joining the points $z_1,$ $z_2.$ Note that the indicated $\inf $ is
calculated and is exactly equal to the right-hand side
in~(\ref{eq3}) (see \cite[theorem~7.2.1, relation~(7.2.5)]{Berd}).

\medskip
If we define the length $l(\gamma)$ of the path $\gamma:[a,
b]\rightarrow {\Bbb D}$ by the equality
\begin{equation}\label{eq30}
l(\gamma):=\sup\limits_{\pi}\sum\limits_{k=0}^nh(\gamma(t_k),
\gamma(t_{k+1}))\,,
\end{equation}
where $h$ is from~(\ref{eq3}), and $\pi=\{a=t_0\leqslant
t_1\leqslant t_2\leqslant\ldots\leqslant t_n=b\},$ then
$l(\gamma)=s_h(\gamma)$ for absolutely continuous paths.

The proof of this fact can be carried out similarly
to~\cite[Theorem~1.3]{Va}, and therefore goes down. Further, for the
rectifiable path $\gamma:[a, b]\rightarrow {\Bbb D}$ there is a
unique normal representation $\gamma_0:[0, l(\gamma)]\rightarrow
{\Bbb D}$ such that $\gamma_0(s(t))=\gamma(t)$ for any $t\in [a,
b],$ where $s(t)$ is the length of the path $\gamma$ on the segment
$[a, t]$ (see, for example, \cite[Sec.~7.1]{He} or
\cite[Sec.~2]{Va}). If the path is only locally rectifiable, then
$l(\gamma)=\infty$ and, accordingly, $\gamma_0:[0,
\infty)\rightarrow {\Bbb D},$ $\gamma_0(s(t))=\gamma(t)$ for any
$t\in [a, b].$ Let $\rho:{\Bbb D}\rightarrow {\Bbb R}$ be a
nonnegative Borel function. Then the integral from $\rho $ over the
(locally rectifiable) path $\gamma$ can be defined by equality
\begin{equation}\label{eq31}
\int\limits_{\gamma}\rho(x)\,ds_h(x)=\int\limits_{0}^{l(\gamma)}\rho(\gamma_0(s))\,ds\,.
\end{equation}
Observe that, the integral in~(\ref{eq31}) coincides with
$$\int\limits_{\gamma}\rho(z)\frac{2\,|dz|}{1-|z|^2}=
\int\limits_{a}^{b}\frac{2\rho(\gamma(t))\,|\gamma^{\,\prime}(t)|}{1-|\gamma(t)|^2}\,dt$$
for absolutely continuous path $\gamma$
(see~\cite[Corollary~2.1]{Sev$_2$}).

\medskip
It should be noted that the normal representation of $\gamma_0(s)$
by an arbitrary locally rectifiable path $\gamma(t)$ is trivial
locally absolutely continuous with respect to its natural parameter
$s\in [0, l(\gamma)]$ ($s\in [0, l(\gamma)),$ if
$l(\gamma)=\infty$). In particular,
$$\int\limits_{0}^{l(\gamma)}\rho(\gamma_0(s))\,ds=\int\limits_{0}^{l(\gamma)}\frac{2\rho(\gamma_0(s))
\,|\gamma_0^{\,\prime}(s)|}{1-|\gamma_0(s)|^2}\,ds\,.$$
In view of the above, we will not distinguish between the length of
the path (integral over it) in (\ref{eq30})--(\ref{eq31}) and
in~(\ref{eq1}), respectively.
\end{remark}

\medskip
In order to simplify research, we introduce into consideration the
so-called {\it fundamental set} $F.$ We define it as a subset of
${\Bbb D},$ containing one and only one point of the orbit $z\in
G_{z_0}$ (see~\cite[\S\,9.1, Ch.~9]{Berd}). A {\it fundamental
domain} $D_0$ is a domain in ${\Bbb D}$ with the property $D_0
\subset F \subset \overline{D_0}$ such that $h(\partial D_0)=0$ (see
ibid). The existence of fundamental sets and fundamental domains is
justified by the presence of their examples, the most important of
which is {\it Dirichlet polygon},
\begin{equation}\label{eq4}
D_{\zeta}=\bigcap\limits_{g\in G, g\ne I}H_g(\zeta)\,,
\end{equation}
where $H_g(\zeta)=\{z\in {\Bbb D}: h(z, \zeta)<h(z, g(\zeta))\}$
(see \cite[relation~(2.6)]{RV}). Let $\pi$ be the natural projection
of ${\Bbb D}$ onto ${\Bbb D}/ G,$ then $\pi$ is an analytic function
conformal on $D_0$ (see also~\cite[Proposition~9.2.2]{Berd} and
comments after~(2.11) in~\cite{RV}). Note, in addition, that there
is a one-to-one correspondence between the points of $F$ and ${\Bbb
D}/G$. For $z_1, z_2 \in F$ we put
\begin{equation}\label{eq5}
d(z_1, z_2):=\widetilde{h}(\pi(z_1), \pi(z_2))\,,
\end{equation}
where $\widetilde{h}$ is defined in~(\ref{eq2}). Observe that
$d(z_1, z_2)\leqslant h(z_1, z_2)$ and, moreover, for any
compactum~$A\subset {\Bbb D}$ there is $\delta=\delta(A)>0$ such
that
\begin{equation}\label{eq34}
d(z_1, z_2)=h(z_1, z_2), \quad \forall\,\, z_1, z_2\in A: h(z_1,
z_2)<\delta\,,
\end{equation}
see e.g.~\cite[Lemma~2.3]{Sev$_2$}.

\medskip
Note that the metric space $(F, d)$ is homeomorphic to $({\Bbb D}/G,
\widetilde{h}).$ Define the elements of length and volume on $(F,
d)$ according to the relations~(\ref{eq1}), in addition, we also
define the elements $ds_{\widetilde{h}}$ and $d\widetilde{h}$ of
length and area on the surface ${\Bbb S},$ respectively as the
elements $ds_h$ and $dh$ in their respective local coordinates.
These local coordinates can be, in particular, fundamental domains
$D_0$ in ${\Bbb D}.$ Due to $d(z_1, z_2)\leqslant h (z_1, z_2) $
and~(\ref{eq34}), the lengths of the paths in the metrics $h$ and
$d$ of the domain $D_0$ coincide.

\medskip
Here and below, $B(z_0, r)$ and $S(z_0, r)$ denote the Euclidean
disk and a circle on a plane centered at the point $z_0\in {\Bbb C}$
and of a radius $r> 0,$ respectively. Let $p_0\in {\Bbb S}$ and
$z_0\in {\Bbb D}$ be such that $\pi(z_0)=p_0,$ where $\pi$ is the
natural projection of ${\Bbb D}$ onto ${\Bbb D}/G.$ Denote by $D_0$
the Dirichlet polygon centered at the point $z_0,$ and put
$\varphi:=\pi^{\,-1}.$ Note that the mapping $\varphi $ is a
homeomorphism of $({\Bbb S}, \widetilde {h})$ onto $(F, d),$ where
$\widetilde{h}$ is a metric on the surface ${\Bbb S},$ a $d$ is the
above-defined metric on the fundamental set $F,$ $D_0\subset
F\subset \overline{D_0}.$ Without loss of generality, we may also
assume that $z_0=0.$ Indeed, otherwise consider an auxiliary mapping
$g_0(z)=(z-z_0)/(1-z\overline {z_0}),$ having no fixed points inside
the unit disk. Then, if $G$ is a group of linear fractional maps
corresponding to the surface ${\Bbb S},$ then
$G^{\,\prime}=\{g_0\circ g, g\in G \}, $ obviously also corresponds
to ${\Bbb S}$ in the sense that the surface ${\Bbb S}$ is again
conformally equivalent to the factor space ${\Bbb D}/G^{\,\prime}.$
Choose a compact neighborhood $V\subset {\Bbb D} $ of the point
$0\in F \subset{\Bbb D},$ such that $d(x, z)=h(x, z)$ for all $x,
z\in V,$ which is possible due to condition~(\ref {eq34}). In
addition, we choose $V$ so that $V\subset B(0, r_0)$ for some
$0<r_0<1.$ Put $U:=\pi(V).$ In this case, the neighborhood $U$ is
called {\it a normal neighborhood of the point $p_0.$} Note that the
ball $\widetilde{B}(p_0, r)\subset U$ corresponds to the set $B_h(0,
r)\subset{\Bbb C}, $ more precisely,
$$
\widetilde{B}(p_0, r):=\{p\in {\Bbb S}: \widetilde{h}(p,
p_0)<r\}=\{p\in {\Bbb S}: h(\varphi(p), 0)<r\}=$$
\begin{equation}\label{eq2A}
=\left\{p\in {\Bbb S}: |\varphi(p)|
<\frac{e^r-1}{e^r+1}\right\}=\pi\left(B\left(0,
\frac{e^r-1}{e^r+1}\right)\right)\,,
\end{equation}
where $h$ is a hyperbolic metric, see~(\ref{eq3}). Similarly,
$$\widetilde{S}(p_0, r):=\{p\in {\Bbb S}: \widetilde{h}(p,
p_0)=r\}=$$
\begin{equation}\label{eq3A}
=\left\{p\in {\Bbb S}:
|\varphi(p)|=\frac{e^r-1}{e^r+1}\right\}=\pi\left(S\left(0,\frac{e^r-1}{e^r+1}\right)\right)\,.
\end{equation}
Throughout what follows, the normal neighborhood $U$ of the point
$p_0,$ as well as the mapping $\varphi$ and the fundamental set $F,$
we will considered selected and fixed. The following analogue of
Fubini's theorem for Riemannian surfaces holds.

\begin{lemma}\label{lem2}
{\sl Let $U$ be some normal neighborhood of the point $p_0\in {\Bbb
S},$ and let $Q:U\rightarrow [0, \infty]$ be a function measurable
with respect to the measure $\widetilde{h},$ and $d_0:={\rm
dist}\,(p_0,
\partial U):=\inf\limits_{p\in
\partial U}\widetilde{h}(p_0, p).$ Then, for any $0<r_0\leqslant d_0$
\begin{equation}\label{eq6}
\int\limits_{\widetilde{B}(p_0,
r_0)}Q(p)\,d\widetilde{h}(p)=\int\limits_0^{r_0}\int\limits_{\widetilde{S}(p_0,
r)}Q(p)\,ds_{\widetilde{h}}(p)\,dr\,,
\end{equation}
where $d\widetilde{h}(p)$ and $ds_{\widetilde{h}}$ are area and
length elements on ${\Bbb S},$ respectively, see~(\ref{eq1}), and
the disk $\widetilde{B}(p_0, r_0)$ and the circle
$\widetilde{S}(p_0, r)$ are defined in~(\ref{eq2A})
and~(\ref{eq3A}). }
\end{lemma}

\medskip
The assertion of Lemma~\ref{lem2} includes the measurability of
$\psi(r):=\int\limits_{\widetilde{S}(p_0,
r)}Q(p)\,ds_{\widetilde{h}}(p)$ with respect to $r$ on the
right-hand side of the integral in~(\ref{eq6}).

\begin{proof}
According to the definition of a normal neighborhood, the ball
$\widetilde{B}(p_0, r_0)\subset {\Bbb S}$ corresponds to the ball
$B_h(0, r_0)\subset F$ in the hyperbolic metric $h.$ Taking into
account the relation (\ref{eq2A}), $B_h(0, r_0)=B\left(0,
\frac{e^{r_0}-1}{e^{r_0}+1}\right).$ By definition
\begin{equation}\label{eq8}
\int\limits_{\widetilde{B}(p_0,
r_0)}Q(p)\,\,d\widetilde{h}(p)=4\int\limits_{B\left(0,
\frac{e^{r_0}-1}{e^{r_0}+1}\right)}\frac{Q(\pi(z))}{(1-|z|^2)^2}\,\,dm(z)\,.
\end{equation}
We use the classical Fubini theorem on the plane (see, for example,
\cite[Theorem~2.6.2]{Fe} or \cite[Theorem~8.1.III]{Sa}). Using polar
coordinates and applying this theorem, we will have that
$$I=4\int\limits_0^{\frac{e^{r_0}-1}{e^{r_0}+1}}\int\limits_{S(0,
r)} \frac{Q(\pi(z))}{(1-|z|^2)^2}\,|dz|\,dr=$$
\begin{equation}\label{eq9}=2\int\limits_0^{\frac{e^{r_0}-1}{e^{r_0}+1}}\frac{1}{1-r^2}\int\limits_{S(0,
r)}\frac{2\, Q(\pi(z))}{1-|z|^2}\,|dz|\,dr\,.
\end{equation}
The last relation takes into account that the function
$\int\limits_{S(0, r)}\frac{2 Q(\pi(z))}{1-|z|^2}\,|dz|$ is
measurable by $r$ (which is also part of the statement of the
classical Fubini theorem). Let us make the change
$t=\log\frac{1+r}{1-r}$ in the last integral in accordance
with~\cite[Theorem~3.2.6]{Fe}. Since $dt=\frac{2 \, dr}{1-r^2},$ we
get:
$$2\int\limits_0^{\frac{e^{r_0}-1}{e^{r_0}+1}}\frac{1}{1-r^2}\int\limits_{S(0,
r)}\frac{2
Q(\pi(z))}{1-|z|^2}\,|dz|\,dr=\int\limits_0^{r_0}\int\limits_{S\left(0,
\frac{e^r-1}{e^r+1}\right)}\frac{2
Q(\pi(z))}{1-|z|^2}\,|dz|\,dr\,=$$
\begin{equation}\label{eq10}
=\int\limits_0^{r_0}\int\limits_{S_h(0,
r)}Q(\pi(z))\,ds_h(z)\,dr=\int\limits_0^{r_0}\int\limits_{\widetilde{S}(p_0,
r)}Q(p)ds_{\widetilde{h}}(p)\,dr=\int\limits_0^{r_0}\psi(r)\, dr\,.
\end{equation}
In particular, by~\cite[Theorem~3.2.6]{Fe} the function $\psi(r)$ is
measurable by~$r.$ Now combining~(\ref{eq8}), (\ref{eq9}) and
(\ref{eq10}), we obtain that
$$\int\limits_{\widetilde{B}(p_0, r_0)}Q(p)\,d\widetilde{h}(p)=\int\limits_0^{r_0}\int\limits_{\widetilde{S}(p_0,
r)}Q(p)\,ds_{\widetilde{h}}(p)\,dr\,,$$
as required to prove.~$\Box$
\end{proof}

Before proceeding directly to the study of mappings on Riemannian
surfaces (including mappings with~(\ref{eq1A})), we formulate the
following statement, which relates the concept of ''almost all''
with respect to the modulus of families of paths and the Lebesgue
sense (its proof is similar to~\cite[Lemma~4.1]{IS}).

\begin{lemma}\label{lem8.2.11}
{\sl\, Let $D$ be a domain of the Riemannian surface ${\Bbb S},$
$p_0 \in \overline {D}$ and let $U$ be some normal neighborhood of
the point $p_0.$ If some property $P$ holds for almost all
intersections $D(p_0, r):=\widetilde{S}(p_0, r)\cap D$ of circles
$\widetilde{S}(p_0, r)$ with a domain $D,$ lying in $U,$ where
''almost all'' is understood in the sense of the modulus of families
of paths and the set
$$E=\{r\in {\Bbb R}: P\,\,\,\,\text{имеет\,\, место для}\,\,\,\,
\widetilde{S}(p_0, r)\cap D\}$$
is Lebesgue measurable, then $P$ also holds for almost all $D(p_0,
r)$ in $U$ with respect to the linear Lebesgue measure by a
parameter $r\in {\Bbb R}.$ Conversely, if $P$ holds for almost all
$D(p_0, r):=\widetilde {S}(p_0, r)\cap D $ with respect to the
linear Lebesgue measure in $r\in {\Bbb R},$ then $P$ also holds for
almost all $D(p_0, r):=\widetilde{S}(p_0, r)\cap D$ in the sense of
a modulus.}
 \end{lemma}

The proof of the following statement is similar
to~\cite[Theorem~9.2]{MRSY} (see also~\cite[Lemma~4.2]{IS}), and is
therefore omitted.

\begin{lemma}\label{lem4A}
{\sl\, Let $D$ and $D_{\,*}$ be domains in ${\Bbb S}$ and ${\Bbb
S}_*,$ respectively, let $p_0\in \overline{D}$ and let $Q\colon
D\rightarrow (0, \infty)$ be a measurable function. Then $f\colon
D\rightarrow D_{\,*}$ satisfies estimate~(\ref{eq1A}) at the point
$p_0$ if and only if there is $0<d_0<\sup\limits_{p\in
D}\widetilde{h}(p, p_0)$ such that
\begin{equation}\label{eq15}
M(f(\Sigma_{\varepsilon}))\geqslant
\int\limits_{\varepsilon}^{\varepsilon_0}
\frac{dr}{\Vert\,Q\Vert(r)}\quad\forall\
\varepsilon\in(0,\varepsilon_0)\,,\ \varepsilon_0\in(0,d_0)\,,
\end{equation}
where, as above, $\Sigma_{\varepsilon}$ denotes the family of all
intersections of the circles $\widetilde{S}(p_0, r)$ with the domain
$D,$ $r\in(\varepsilon, \varepsilon_0),$
$$
\Vert Q\Vert(r)=\int\limits_{D(p_0,r)}Q(p)\,ds_{\widetilde{h}}(p)$$
is $L_1$-norm of the function $Q$ over the intersection $D\cap
\widetilde{S}(p_0,r)=D(p_0,r)=\{p\in D\,:\, \widetilde{h}(p,
p_0)=r\}.$ }
\end{lemma}

\section{Basic bound for distortion of the modulus of families of
paths}\label{sec3}

We say that a set $A \subset {\Bbb S}$ has {\it Lebesgue measure
zero} if $A$ can be covered by at most countable the number of
normal neighborhoods $U_k \subset {\Bbb S},$ $k=1,2, \ldots, $ such
that $\varphi_k(U_k) \rightarrow {\Bbb D}, $ where $\varphi_k$ is
some homeomorphisms related to each other conformal transformation,
in this case, $m(\varphi_k(U_k \cap A))=0$ for any $k=1,2, \ldots, $
$m$ is the Lebesgue measure in ${\Bbb C}.$ The following statements
are true.

\begin{lemma}\label{pr1}
{\sl\, Suppose that $B_0\subset {\Bbb S}$ has a Lebesgue measure
zero, $p_0\in {\Bbb S},$ $U$ is a normal neighborhood of the point
$p_0,$ $\overline{U}\ne {\Bbb S}$ and $0<\varepsilon_0<{\rm
dist}\,(p_0, \partial U). $ Then
\begin{equation}\label{eq18}
{\mathcal H}^{\,1}(\varphi(B_0\cap S_r))=0
\end{equation}
for almost all circles $S_r:=\widetilde{S}(p_0,r)$ centered at a
point $p_0,$ where $\varphi=\pi^{-1}$ is a homeomorphism of $U$ into
${\Bbb D},$ corresponding to the definition of a normal neighborhood
$U,$ ${\mathcal H}^{\, 1}$ is a 1-dimensional Hausdorff measure in
${\Bbb C},$ and ''almost all'' should be understood with respect to
the parameter $r\in (0, \varepsilon_0).$ }
\end{lemma}

\begin{proof}
Indeed, since the Lebesgue measure is regular, there is a Borel set
$B\subset U$ such that $\varphi(B_0)\subset \varphi(B)$ and
$m(\varphi(B_0))=m(\varphi(B))=0,$ where $m$ is, as usual, the
Lebesgue measure in ${\Bbb C}.$ Let $g$ be the characteristic
function of the set $\varphi (B).$ According
to~\cite[Theorem~3.2.5]{Fe} for $m= 1,$ we have that
\begin{equation}\label{eq19}
\int\limits_{\varphi(\gamma)}g(z)|dz|={\mathcal
H}^{\,1}(\varphi(B\cap |\gamma|))\,,
\end{equation}
where $\gamma:[a, b]\rightarrow U$ is any locally rectifiable path,
$|\gamma|$ us a locus of $\gamma$ in $U,$ and $|dz|$ is an element
of the Euclidean measure. Arguing similarly to the proof of
\cite[Theorem~33.1]{Va}, we put
$$\rho(p)= \left \{\begin{array}{rr}\infty, & p\in B,
\\ 0 ,  &  p \notin B\ .
\end{array} \right.$$
Observe that $\rho$ is a Borel function. Let $\Gamma$ be a family of
all circles $S_r:=\widetilde{S}(p_0,r)$ centered at the point $p_0,$
for which ${\mathcal H}^{\,1}(\varphi(B\cap S_r))>0.$
By~(\ref{eq19}), for any $S_r\in\Gamma$ we obtain that
$$\int\limits_{S_r}\rho(p)\,ds_{\widetilde{h}}(p)=\int\limits_{S_h(0,
r)}\rho(\varphi^{\,-1}(y))\,ds_h(y)=2\int\limits_{S\left(0,
\frac{e^r-1}{e^r+1}\right)}\frac{\rho(\varphi^{\,-1}(y))}{1-|y|^2}|dy|=$$
$$=2\int\limits_{S\left(0,
\frac{e^r-1}{e^r+1}\right)}\frac{g(y)\rho(\varphi^{\,-1}(y))}{1-|y|^2}|dy|=\infty\,.$$
Now $\rho\in{\rm adm}\,\Gamma.$ Thus, $M(\Gamma)\leqslant
\int\limits_{\Bbb S}\rho^2(p)\,d\widetilde{h}(p)=0.$ Let
$\Gamma^{\,*}$ be a family of all circles
$S_r:=\widetilde{S}(p_0,r)$ centered at $p_0$ for which ${\mathcal
H}^{\,1}(\varphi(B_0\cap S_r))>0.$ Observe that
$\Gamma^{\,*}\subset\Gamma,$ whence $M(\Gamma^{\,*})=0.$ Finally,
note that the function $\psi(r):={\mathcal H}^{\,1}(\varphi(B_0\cap
S_r))$ is Lebesgue measurable by the classical Fubini theorem, so
that~(\ref{eq18}) is true for almost all $r\in (0, \varepsilon_0)$
by Lemma~\ref{lem8.2.11}.~$\Box$
\end{proof}

Let $\gamma:[a, b]\rightarrow {\Bbb S}$ be a (locally rectifiable)
path on the Riemannian surface ${\Bbb S}.$ Then we define the
function $l_{\gamma}(t)$ as the length of the path $\gamma|_{[a,
t]},$ $a\leqslant t\leqslant b $ (where ''length'' is understood in
the sense of a Riemannian surface). For a set $B\subset {\Bbb S},$
put
\begin{equation}\label{eq36}
l_{\gamma}(B)={\rm mes}_1\,\{s\in [0, l(\gamma)]: \gamma(s)\in
B\}\,,
\end{equation}
where, as usual, ${\rm mes}_1$ denotes the linear Lebesgue measure
in ${\Bbb R},$ and $l(\gamma)$ is the length of~$\gamma.$ Similarly,
we may define the value $l_{\gamma}(B)$ for the dashed line
$\gamma,$ i.e. when $\gamma:~\bigcup\limits_{i=1}^{\infty}(a_i,
b_i)\rightarrow {\Bbb S},$ where $a_i<b_i$ for any $i\in {\Bbb N}$
and $(a_i, b_i)\cap (a_j, b_j)=\varnothing$ for any $i\ne j.$

\begin{lemma}\label{pr2}
{\sl\, Let $D$ and $D_{\,*}$ be domains in ${\Bbb S}$ and ${\Bbb
S}_*,$ respectively, and let $f:D\rightarrow D_{\,*}$ be a mapping
of the Sobolev class $W^{1, 1}_{\rm loc}.$ Let $p_0\in
\overline{D},$ let $U$ be a normal neighborhood of $p_0,$
$\overline{U}\ne {\Bbb S}$ and $0<\varepsilon_0<{\rm dist}\,(p_0,
\partial U),$ and let $B_0\subset D$ has a Lebesgue measure zero. Then
${\mathcal H}^{\,1}(f(B_0\cap \widetilde{S}(p_0, r)))=0$ for almost
any $r\in (0, \varepsilon_0)$ in local coordinates and, in addition,
\begin{equation}\label{eq32}
l_{f(\widetilde{S}(p_0, r)\cap D)}(f(B_0))=0\,,
\end{equation}
where $l$ is defined in~(\ref{eq36}). }
\end{lemma}

\begin{proof}
Since the mapping $f$ is continuous, the domain $f(D)$ can be
covered by at most a countable number of neighborhoods $V_k,$
$k=1,2, \ldots, $ in such a way, that $V_k$ is conformally
homeomorphic to some neighborhood $W_k \subset {\Bbb D},$
$\overline{W_k}$ is compact in ${\Bbb D}$ and, moreover, $f^{\,-
1}(V_k)=U_k\subset D,$ where $U_k$ is an open set, $\bigcup
\limits_{k=1}^{\infty} U_k=D.$ We may also assume that the length
and area in $V_k$ are calculated in terms of the hyperbolic length
and hyperbolic area in $W_k.$ By what was said above, without loss
of generality, we may assume that $f(D)$ is conformally homeomorphic
to the set $W\subset {\Bbb D},$ whose closure is compact in $ {\Bbb
D}.$ Let the indicated conformal homeomorphism be realized for using
the mapping $\psi: f(D)\rightarrow W, $ and let $\varphi$ be a
homeomorphism of $U$ into ${\Bbb D},$ corresponding to the
definition of normal neighborhood $U,$ and let $\varphi(U)\subset
B(0, r_0),$ where $\overline{B(0, r_0)}$ is a compact set in ${\Bbb
D}.$ Consider a partition of the set $B(0, r_0)$ into a countable
number of pairwise disjoint ring segments
\begin{equation}\label{eq23}
A_m=\{z\in {\Bbb C}: z=Re^{i\alpha}, R\in (r_{m-1}, r_m], \alpha\in
(\psi_{m-1}, \psi_m]\}\,, m\in {\Bbb N}\,.
\end{equation}
Let $h_m$ be an auxiliary quasiisometry that maps $A_m$ onto a
rectangle~$B_m$ such that arcs of circles centered at zero are
mapped to line segments, see Figure~\ref{fig1}.
\begin{figure}[h]
\centerline{\includegraphics[scale=0.5]{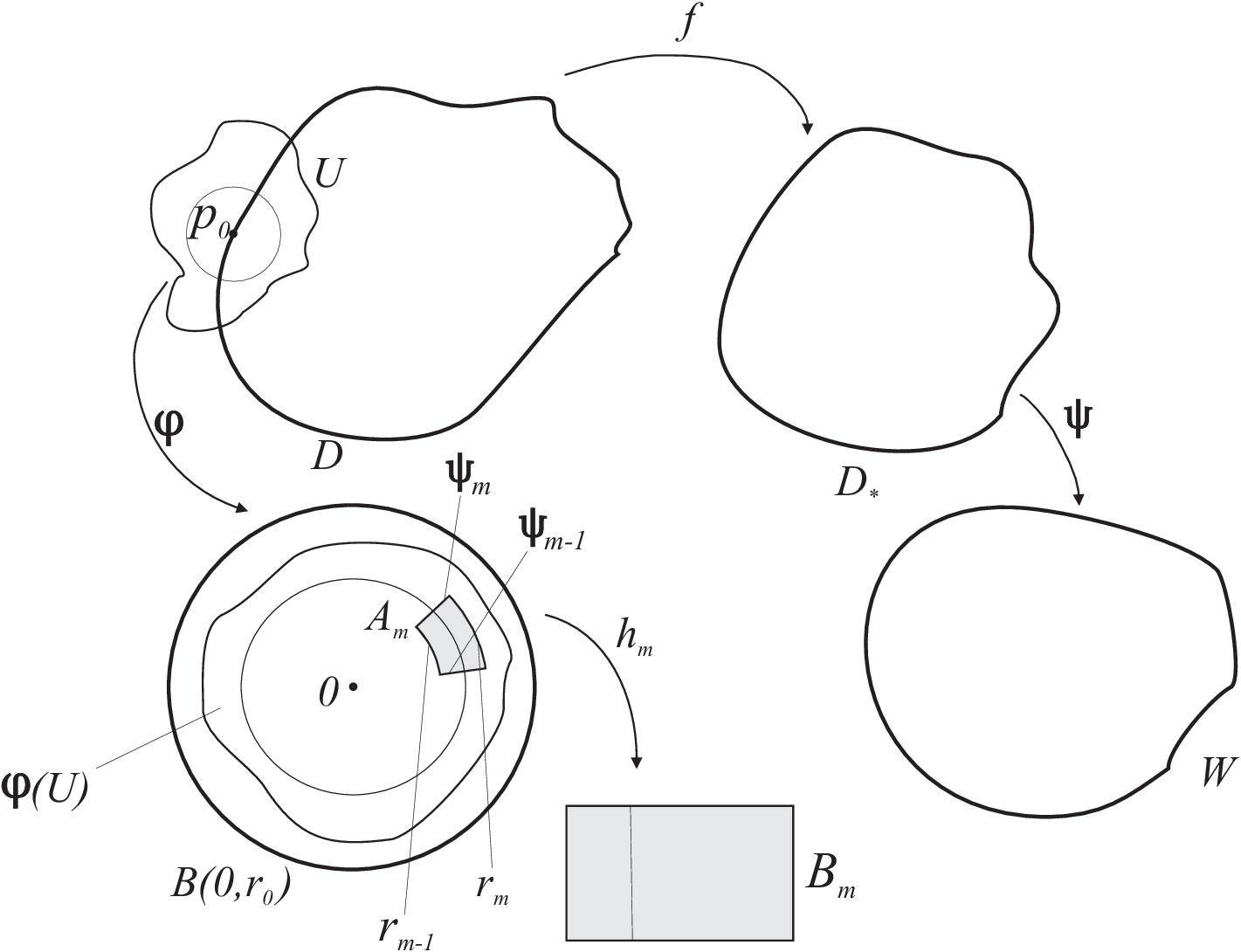}} \caption{To
the proof of Lemma~\ref{pr2}.}\label{fig1}
\end{figure}
More precisely, put $h_m(\omega)=\log\omega,$ $\omega\in A_m$ and
$\widetilde{A_m}:=h_m(A_m\cap D).$ Then, for each $m\in {\Bbb N},$
consider the mapping
$$g_m:=\psi\circ
f\circ\varphi^{\,-1}\circ h_m^{\,-1}\,,\quad
g_m:\widetilde{A_m}\rightarrow {\Bbb C}\,.$$
Observe that $g_m\in W^{1, 1}_{\rm loc}(\widetilde{A_m})$ (see
\cite[Section~1.1.7]{Ma}), whence, in particular, $g_m\in ACL$ (see
\cite[Theorems~1 and 2, Section~1.1.3, $\S\,$1.1, Section.~I]{Ma}).
Set, as above, $S_r:=\widetilde{S}(p_0, r).$ By Lemma~\ref{pr1} and
in view of the smoothness of the mapping $h_m,$ we obtain that
$${\mathcal H}^{\,1}(\varphi(B_0\cap S_r)\cap A_m)={\mathcal
H}^{\,1}(h_m(\varphi(B_0\cap S_r)\cap A_m))=0$$
for any $r\in [0, \varepsilon_0]\setminus A_0,$ where ${\rm
mes}_1A_0=0.$ Set $D_r:=D\cap S_r.$ Then also for any~$r\in [0,
\varepsilon_0]\setminus A_0$
\begin{equation}\label{eq20}
{\mathcal H}^{\,1}(\psi(f(B_0\cap D_r\cap
\varphi^{-1}(A_m))))={\mathcal H}^{\,1}(g_m(h_m(\varphi(B_0\cap
D_r)\cap A_m)))=0\,,
\end{equation}
since the absolute continuity of the map $g_m$ on a fixed interval
implies the $N$ -property with respect to the linear Lebesgue
measure (see~\cite[Section~2.10.13]{Fe}). Observe that $U\subset
\bigcup\limits_{m=1}^{\infty}\varphi^{-1}(A_m),$ so
from~(\ref{eq20}), in view of the countable semi-additivity of the
Hausdorff measure,
\begin{equation}\label{eq21}
{\mathcal H}^{\,1}(\psi(f(B_0\cap D_r)))=0\,, \quad r\in [0,
\varepsilon_0]\setminus A_0\,.
\end{equation}
Let $\gamma_i$ be an arbitrary arc of the dashed line
$\psi(f(D_r)).$ Then we parametrize $\gamma_i:[0,
l(\gamma_i)]\rightarrow {\Bbb D},$ $\gamma_i=\gamma_i(s),$ where $s$
is a natural parameter on $\gamma_i$ in the sense of the Euclidean
length. Setting $m=1$ in~\cite[Theorem~3.2.5]{Fe}, by~(\ref {eq21})
we obtain that the set $B:=\{s\in [0, l(\gamma_i)]: \gamma_i(s)\in
\psi(f(B_0))\}$ has a linear measure zero. Let
$\chi_{\psi(f(B_0))}(z)$ be a characteristic function of the set
$\psi(f(B_0)).$ Taking into account the Remark~\ref {rem1}, we
obtain that
$$l_{f(\widetilde{S}(p_0,
r)\cap D )}(f(B_0))=\sum\limits_{i=1}^{\infty}
2\int\limits_0^{l(\gamma_i)}\frac{\chi_{\psi(f(B_0))}(\gamma_i(s))\,
ds}{1-|\gamma_i(s)|^2}=0$$
for almost any $r\in (0, \varepsilon_0),$ which was required to be
established.~$\Box$
\end{proof}

\begin{lemma}\label{pr3}
{\sl\, Let $D$ and $D_{\,*}$ be domains in Riemannian surfaces
${\Bbb S}$ and ${\Bbb S}_*,$ respectively, and let $f:D\rightarrow
D_{\,*}$ be an open discrete $W^{1, 1}_{\rm loc}$-mapping  with a
finite distortion. Let $p_0\in \overline{D},$ let $U$ be a normal
neighborhood of the point $p_0,$ let $0<\varepsilon_0<{\rm
dist}\,(p_0, \partial U),$ and let $B_*\subset D $ be the set of
such points of  $D,$ in which the mapping $f$ is differentiable (in
local coordinates), however, $J_f(p)=0.$ Then
\begin{equation}\label{eq33}
l_{f(\widetilde{S}(p_0, r)\cap D)}(f(B_*))=0
\end{equation}
for almost all $r \in (0, \varepsilon_0),$
where the function $l$ is defined in~(\ref{eq36}).}
\end{lemma}

\begin{proof}
Observe that $f$ is differentiable almost everywhere in $D$ in local
coordinates (see the remarks made in the introduction in~\cite{RV}).
In particular, the set $U$ can be split into a countable number sets
$B_k,$ $k=0,1,2, \ldots, $ such that $f|_{B_k}$ is a bilipschitz
homeomorphism for $k=1,2, \ldots, $ and $B_0$ has a measure zero
(see~\cite[items~3.2.2, 3.1.4 and 3.1.8]{Fe}). Let, as before, $S_r:
=\widetilde{S}(p_0, r)$ and $D_r:=S_r \cap D.$ By Lemma~\ref{pr2}
${\mathcal H}^{\,1}(f(B_0\cap D_r))=0$ for almost of any $r\in (0,
\varepsilon_0)$ in local coordinates, therefore, a 1-dimensional
change of variables holds for almost all $r\in (0, \varepsilon_0),$
(see~\cite[Theorem~3.2.5 ]{Fe}).

\medskip
Repeating the reasoning given in the proof of Lemma~\ref{pr2} and
using the notation of this proposition, we conclude that the mapping
$h_m$ maps $\varphi(S_r)\cap A_m $ to some part of the segment $I(m,
R)=\{z \in {\Bbb C}: z=\log R + it, t \in (\psi_{m-1}, \psi_m), R =
(e^r-1)/(e^r+ 1)\}. $ Since $f$ has a finite distortion,
$g_m^{\,\prime}(\log R+ it)=0$ for all $t\in(\psi_{m-1}, \psi_m)$
such that $\varphi^{\,- 1}(h_m^{\,- 1}(\log R+ it)) \in B_*. $ Then,
by virtue of~\cite[Theorem~3.2.5]{Fe} and in view of the above
remarks
$${\mathcal H}^{\,1}(\psi(f(B_*\cap D_r\cap
\varphi^{-1}(A_m))))={\mathcal H}^{\,1}(g_m(h_m(\varphi(B_*\cap
D_r)\cap A_m)))\leqslant$$
$$\leqslant\int\limits_{g_m(h_m(\varphi(B_*\cap
D_r)\cap A_m))}N(y, g_m, h_m(\varphi(B_*\cap D_r)\cap
A_m))\,d{\mathcal H}^{1}y=$$$$=
\int\limits_{\psi_{m-1}}^{\psi_m}\chi_{h_m(\varphi(B_*)\cap
A_m)}(\log R+it)|g_m^{\,\prime}(\log R+it)|dt=0\,,$$
where $\chi_{h_m(\varphi(B_*)\cap A_m)}$ is a characteristic
function of the set $h_m(\varphi(B_*)\cap A_m).$ Semiadditivity with
respect to $m$ of the one-dimensional Hausdorff measure in the last
chain of equalities gives us ${\mathcal H}^{\,1}(\psi(f(B_*\cap
D_r)))=0$ for almost all $r\in (0, \varepsilon_0). $ Let $\gamma_i$
be an arbitrary dashed arc line $\psi(f(D_r)).$ Parametrize
$\gamma_i$ as $\gamma_i:[0, l(\gamma_i)]\rightarrow {\Bbb D},$
$\gamma_i=\gamma_i(s),$ where $s\in [0, l (\gamma_i)]$ is a natural
parameter. Setting $m=1$ in \cite[Theorem ~ 3.2.5]{Fe}, we obtain
that the set $B_i:=\{s\in [0, l(\gamma_i)]: \gamma_i(s)\in
\psi(f(B_*))\}$ has a linear measure zero. Let
$\chi_{\psi(f(B_*))}(z)$ be the characteristic function of the set
$\psi(f(B_*)).$ Taking into account the Remark~\ref{rem1}, we obtain
that
$$l_{f(\widetilde{S}(p_0,
r)\cap
D)}(f(B_*))=\sum\limits_{i=1}^{\infty}2\int\limits_0^{l(\gamma_i)}\frac{\chi_{\psi(f(B_*))}(\gamma_i(s))\,
ds}{1-|\gamma_i(s)|^2}=0$$
for almost any $r\in (0, \varepsilon_0),$ which was required to be
established.~$\Box$
\end{proof}

\begin{figure}[h]
\centerline{\includegraphics[scale=0.5]{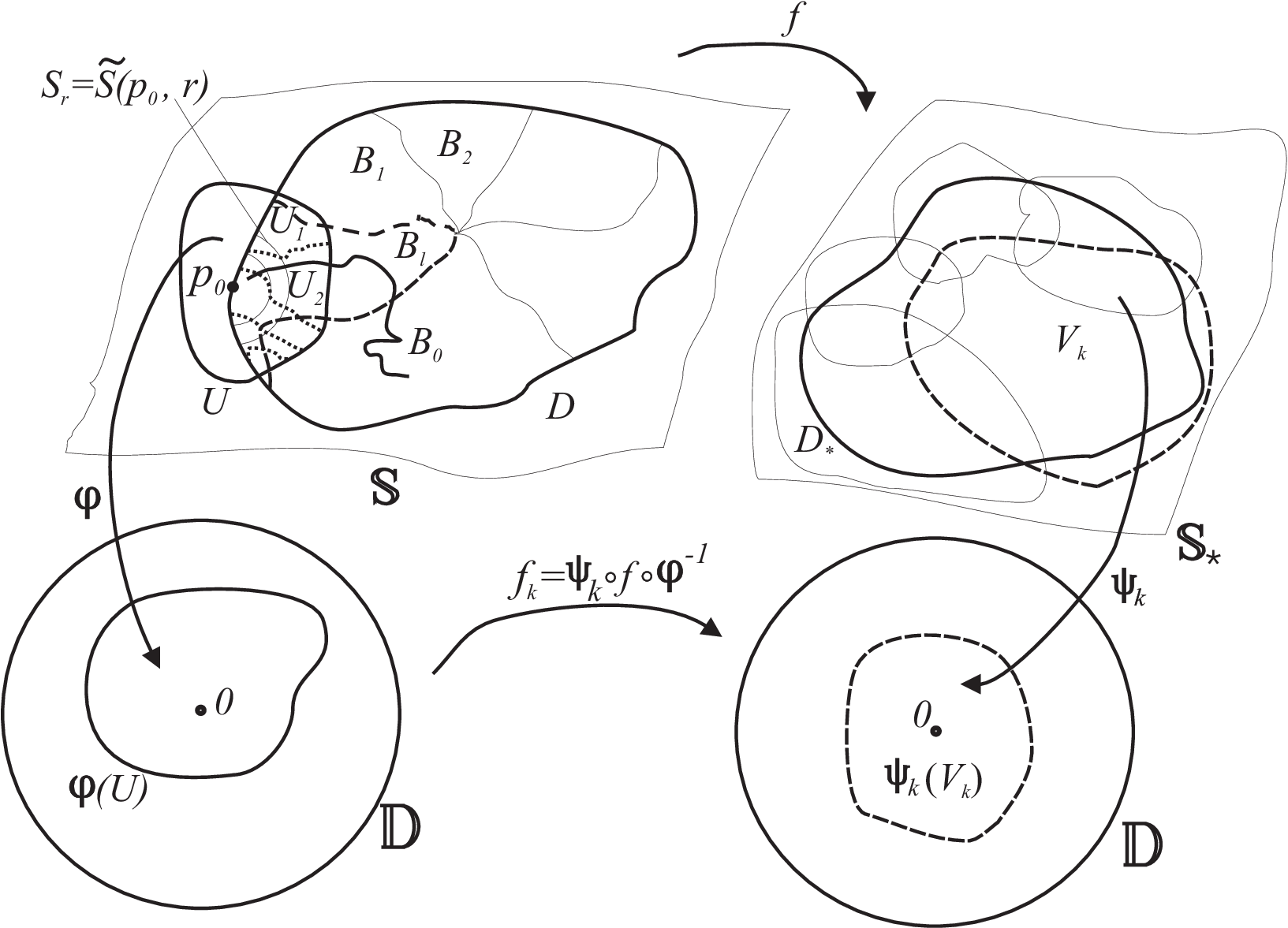}} \caption{To
the proof of Theorem~\ref{thOS4.2}}\label{fig4}
\end{figure}

{\it Proof of Theorem~\ref{thOS4.2}}. Since $f$ is open, the mapping
$f$ is differentiable almost everywhere in $D$ local coordinates
(see the remarks made in the introduction to~\cite {RV}; see
also~\cite[Theorem~III.3.1]{LV}). Let $ B $ be the Borel set of all
points $p\in D,$ where $f$ has a total differential $f^{\,\prime}(p)
$ and $J_f(p)\ne 0$ in local coordinates. Note that $B$ may be
represented as at most countable unions of Borel sets $B_l,$ $l=1,2,
\ldots\,,$ such that $f_l=f|_{B_l}$ are bilipschitz homeomorphisms
(see~\cite[Sections~3.2.2, 3.1.4 and 3.1.8]{Fe}). See
Figure~\ref{fig4} for illustrations.
Without loss of generality, we may assume that the sets $B_l$ are
pairwise disjoint. We also denote by $B_*$ the set of all points
$p\in D,$ where $f$ has a total differential and
$f^{\,\prime}(p)=0.$

Since $f$ is a finite distortion, $f^{\,\prime}(p)=0$ for almost all
points $p,$ where $J_f(p)=0.$ Thus, by construction, the set
$B_0:=D\setminus \left(B\bigcup B_*\right)$ has $\widetilde
{h}$-measure zero.  Let $U$ be a normal neighborhood of the point
$p_0$ and $\varphi:U\rightarrow {\Bbb D}$ be a mapping corresponding
to this normal neighborhood. We may assume that $\varphi(U) \subset
B(0, r_0),$ $0<r_0<1.$ Since $\overline{D_ *}$ is compact in ${\Bbb
S}_*,$ we can cover $\overline{D_*}$ by a finite number of
neighborhoods $V_k,$ $k=1,2,\ldots, m, $ such that $\psi_k: V_k
\rightarrow B(0, R_k),$ $0<R_k <1, $ and $\psi_k $ are conformal
mappings. Let $R_0:=\max\limits_{1\leqslant k \leqslant m} R_k. $
Since the mapping $f$ is continuous, the sets
$U^{\,\prime}_k:=f^{\,-1}(V_k\cap D_{\,*})\cap U$ are open in $U$
and the mapping
$$f_k:=\psi_k\circ f\circ\varphi^{\,-1}$$
is a mapping from $\varphi(U^{\,\prime}_k)\subset B(0, r_0)$ into
$\psi_k(V_k)\subset B(0, R_0).$

\medskip
Set $U_1=U^{\,\prime}_1,$ $U_2=U^{\,\prime}_2\setminus
U^{\,\prime}_1,$ $U_3=U^{\,\prime}_3\setminus (U^{\,\prime}_1\cup
U^{\,\prime}_2),$ $\ldots, U_m=U^{\,\prime}_m\setminus
(U^{\,\prime}_1\cup U^{\,\prime}_2\ldots U^{\,\prime}_{m-1}).$
Observe that, by the definition, $U_m\subset U^{\,\prime}_m$ for
$m\geqslant 1$ and $U_s\cap U_k=\varnothing$ for $s\ne k.$ Let
$\Gamma$ be a family $D_r:=D\cap S_r$ of all intersections of
circles $S_r=\widetilde{S}(p_0, r),$ $r\in(\varepsilon, r_0),$ with
$D.$ We fix an admissible function $\rho_*\in{\rm adm}\,f(\Gamma),$
$\rho_*\equiv 0$ outside $f(D)$, and put $\rho \equiv 0$ outside $U$
and on $B_0,$ and
$$\rho(p)\colon=\rho_*(f(p))\Vert f_k^{\,\prime}
(\varphi(p))\Vert \qquad\text{при}\ p\in U_k\setminus B_0\,,$$
where the matrix norm of the derivative $\Vert g^{\,\prime}(z)\Vert$
of a given function $g:{\Bbb D}\rightarrow {\Bbb C},$ ${\Bbb
D}\subset {\Bbb C},$ as usual, is defined as $\Vert
g^{\,\prime}(z)\Vert=|g_z|+|g_{\overline{z}}|,$ $g_z=(g_x-ig_y)/2,$
$g_{\overline{z}}=(g_x+ig_y)/2,$ $z=x+iy\in {\Bbb C}.$ Observe that
$$D_r:=D\cap
S_r=\left(\bigcup\limits_{1\leqslant k\leqslant m\atop 1\leqslant
l<\infty} D^{\,r}_{kl}\right)\bigcup \left(\bigcup\limits_{k=1}^m
S_r\cap U_k\cap B_*\right)\bigcup \left(\bigcup\limits_{k=1}^m
S_r\cap U_k\cap B_0\right)\,,$$
where $D^{\,r}_{kl}=S_r\cap U_k\cap B_l.$ By Lemmas~\ref{pr2}
and~\ref{pr3} $l_{f(D_r)}(f(U_k\cap B_0))=0$ and
$l_{f(D_r)}(f(U_k\cap B_*))=0$ for any $1\leqslant k\leqslant m$ and
almost any $r\in (0, r_0).$ Thus,
$$1\leqslant\int\limits_{D^{\,*}_r}\rho_*(p_*)\,ds_{\widetilde{h_*}}(p_*)\leqslant
2\sum\limits_{k=1}^m
\sum\limits_{l=1}^{\infty}\,\,\,\int\limits_{\psi_k(f(D^{\,r}_{kl}))}
\frac{\rho_*(\psi_k^{\,-1}(z))}{1-|z|^2}\,d{\mathcal
H}^{\,1}(z)\leqslant$$
\begin{equation}\label{eq27}\leqslant \frac{2}{1-R_0^2}\sum\limits_{k=1}^m
\sum\limits_{l=1}^{\infty}\,\,\,\int\limits_{\psi_k(f(D^{\,r}_{kl}))}
\rho_*(\psi_k^{\,-1}(z))\,d{\mathcal H}^{\,1}(z)
\end{equation}
for almost any $r\in (0, r_0),$ where ${\mathcal H}^{\,1}$ denotes
1-dimensional Hausdorff measure. On the other hand, arguing at each
set~$D^{\,r}_{kl}$ separately and using \cite[item~1.7.6,
theorem~2.10.43 and theorem~3.2.6]{Fe}, we obtain that
$$\int\limits_{D^{\,r}_{kl}}\rho(p)\,ds_{\widetilde{h}}(p)=
2\int\limits_{\varphi(D^{\,r}_{kl})}\frac{\rho_*(\psi^{\,-1}_k(f_k(y))\Vert
f_k^{\,\prime}(y)\Vert}{1-|y|^2}\,\,|dy|\geqslant$$
\begin{equation}\label{eq29}
2\int\limits_{\varphi(D^{\,r}_{kl})}\rho_*(\psi^{\,-1}_k(f_k(y))\Vert
f_k^{\,\prime}(y)\Vert\,\,|dy|\geqslant
2\int\limits_{\psi_k(f(D^{\,r}_{kl}))}\rho_*(\psi_k^{\,-1}(z))\,d{\mathcal
H}^{\,1}(z)\,.
\end{equation}
Summing~(\ref{eq29}) over all $1\leqslant k\leqslant m$ and
$1\leqslant l<\infty,$ and considering~(\ref{eq27}) and
Lemma~\ref{lem8.2.11}, we conclude
that~$\rho/(1-R_0^2)\in{\rm{ext\,adm}}\,\Gamma.$

\medskip
Using the change of variables on each $B_l,$ $l=1,2,\ldots $ (see,
e.g. \cite[Theorem~3.2.5]{Fe}), countable additivity of the Lebesgue
integral, and also taking into account~(\ref{eq1}), we obtain the
estimate
$$\frac{1}{(1-R_0^2)^2}\int\limits_{D}\frac{\rho^2(p)}{K_{\mu}(p)}\,d\widetilde{h}(p)=$$$$=
\frac{4}{(1-R_0^2)^2}\sum\limits_{k=1}^m
\sum\limits_{l=1}^{\infty}\int\limits_{\varphi(U_k\cap
B_l)}\frac{\rho^2_*((f\circ \varphi^{-1})(z))\Vert f_k^{\,\prime}
(z)\Vert^2}{(1-|z|^2)^2K_{\mu}(\varphi^{\,-1}(z))}\,dm(z)\leqslant$$
$$\leqslant\frac{4}{(1-r^2_0)^2(1-R_0^2)^2}\sum\limits_{k=1}^m
\sum\limits_{l=0}^{\infty}\int\limits_{{\Bbb
S}_*}\rho^2_*(\psi_k^{\,-1}(y))N(y, f_k, \varphi(U_k\cap
B_l))\,dm(y)\leqslant $$$$ \leqslant\frac{4}{(1-r^2_0)^2(1-R_0^2)^2}
\sum\limits_{k=1}^m\int\limits_{{\Bbb
S}_*}\rho^2_*(\psi_k^{\,-1}(y))N(y, f_k,
\varphi(U_k))\,dm(y)\leqslant$$$$\leqslant
\frac{4}{(1-r^2_0)^2(1-R_0^2)^2}
\sum\limits_{k=1}^m\int\limits_{{\Bbb
S}_*}\frac{\rho^2_*(\psi_k^{\,-1}(y))N(\psi^{\,-1}_k(y), f,
U_k)}{(1-|y|^2)^2}\,dm(y)=$$
$$=\frac{1}{(1-r^2_0)^2(1-R_0^2)^2}
\sum\limits_{k=1}^m\int\limits_{{\Bbb S}_*}\rho^2_*(p_*)N(p_*, f,
U_k)\,d\widetilde{h_*}(p_*)\leqslant$$
$$\leqslant\frac{N(f,
D)}{(1-r^2_0)^2(1-R_0^2)^2} \int\limits_{{\Bbb
S}_*}\rho^2_*(p_*)\,d\widetilde{h_*}(p_*)\,.$$
To complete the proof, one should put
$c:=\frac{1}{(1-r^2_0)^2(1-R_0^2)^2}.~\Box$

\section{Boundary extension of mappings with lower modulus distortion}\label{sec4}

So, we have established the main modulus inequality for the Sobolev
classes, with which we will work further. Now let us talk on the
boundary extension of the Sobolev classes, for which we consider an
auxiliary class of mappings with the condition~(\ref{eq1A}).

We recall the following definitions. A mapping $f:D\rightarrow
D_{\,*},$ $f(D)=D_{\,*},$ is called {\it boundary preserving}, if
$C(f,
\partial D)\subset
\partial D_{\,*},$ where, as usual, $$C(f, \partial D)=\{p_*\in{\Bbb
S}_*: \exists\,\, p_k\in D, p\in \partial D: p_k\rightarrow p,
f(p_k)\rightarrow p_*, k\rightarrow\infty\}\,.$$
The following statement is established in~\cite[Theorem~3.3]{Vu$_1$}
for the case of the space ${\Bbb R}^n.$ In our case, its validity
directly follows from~\cite[Proposition~2.1]{Sev$_1$}.

\begin{proposition}\label{pr4}
{\,\sl Let $D$ and $D_{\,*}$ be domains in Riemannian surfaces
${\Bbb S}$ and ${\Bbb S}_*,$ respectively. Let $f:D\rightarrow
D_{\,*}$ be open discrete and closed mapping in $D.$ Then $f$ is
boundary preserving. }
\end{proposition}

A Borel function $\rho:{\Bbb D}\rightarrow [0, \infty]$ is called
{\it admissible} for the family $\Gamma$ of paths $\gamma$ in the
sense of hyperbolic length, write $\rho\in {\rm adm}_h\,\Gamma,$ if
$\int\limits_{\gamma}\rho(z)\,ds_h(z)\geqslant 1$ for any path
$\gamma\in \Gamma,$ where $ds_h(z)$ is an element of length
corresponding to~(\ref{eq1}). {\it The modulus} of the family
$\Gamma$ in the sense of a hyperbolic measure is the quantity
$M_h(\Gamma):=\inf\limits_{\rho\in {\rm
adm}_h\,\Gamma}\int\limits_{\Bbb D}\rho^2(z)\,dh(z),$ where $dh(z)$
is an element of the hyperbolic area. The following result holds,
see~\cite[Remark~5.2]{Sev$_1$}.

\begin{proposition}\label{pr6}
{\sl\, Let $\Gamma$ be a family of paths in ${\Bbb D}.$ Now
$$M_h(\Gamma)=M_e(\Gamma)\,.$$
}
\end{proposition}
An analogue of the following statement is established for the space
${\Bbb R}^n$ in V\"{a}is\"{a}l\"{a}'s
monograph~\cite[Theorem~7.5]{Va}.

\begin{proposition}\label{pr5}{\sl\,
Let ${\Bbb S}$ be a Riemannian surface and $p_0\in {\Bbb S}.$ Let
$U$ be a neighborhood of the point $p_0$ such that $\widetilde{h}(p,
p_0)=h(\varphi(p), \varphi(p_0)),$  where $\varphi=\pi^{\,-1}$ and
$\pi$ is the natural projection of the fundamental Dirichlet polygon
$D_0$ with center at the point $\varphi(p_0)$ on ${\Bbb S}.$ Let
$0<r_1<r_2<{\rm dist}\,(p_0, \partial U),$
$\widetilde{S_i}=\widetilde{S}(p_0, r_i),$ $i=1,2,$
$\widetilde{A}(p_0, r_1, r_2)=\{p\in {\Bbb S}:
r_1<\widetilde{h}(p_0, p)<r_2\}.$ If $\Gamma=\Gamma(\widetilde{S_1},
\widetilde{S_2}, \widetilde{A})$  is family of paths joining
$\widetilde{S_1}$ and $\widetilde{S_2}$ in $\widetilde{A},$ then
$$M(\Gamma)=\frac{2\pi}{L(r_1, r_2)}\,,$$
where $L(r_1,
r_2)=\left(\frac{e^{r_2}-1}{e^{r_2}+1}\right):\left(\frac{e^{r_1}-1}{e^{r_1}+1}\right).$
 }
\end{proposition}

\begin{proof}
By the definition of the mapping $\varphi$ and the neighborhood $U,$
$M(\Gamma)=M_h(\Gamma^{\,*}),$ where
$\Gamma^{\,*}=\Gamma(S_h(\varphi(p_0), r_1), S_h(\varphi(p_0), r_2),
A_h),$ $A_h:=\{z\in {\Bbb D}: r_1<h(z, \varphi(p_0))<r_2 \}.$ By
Proposition~\ref{pr6} $M(\Gamma)=M_e(\Gamma^{\,*}).$ The required
conclusion now follows from~\cite[Corollary~5.18]{Vu$_2$}.~$\Box$
\end{proof}

\medskip
Let $\Omega$ be a domain in ${\Bbb C},$ or a domain in ${\Bbb S}.$
According to~\cite[item~3]{Na}, the connected set $E\subset\Omega$
is called {\it cut} if $E$ is closed in $\Omega,$ $\overline{E}\cap
\partial \Omega \ne \varnothing $ and $\Omega \setminus E $
consists of two components, the boundary of each of which intersects
$\partial \Omega. $ A sequence of cuts $E_1, E_2, \ldots, E_k,
\ldots $ is called a {\it chain} if $E_k$ separates $E_{k-1}$ from
$E_{k+1}$ in $ \Omega, $ that is, $E_{k-1}$ and $E_{k+1}$ belong to
different components of $\Omega \setminus E_k.$ It follows from the
above definitions that one of the subdomains $\Omega \setminus E_k $
contains all $E_m$ for $m> k.$ This subdomain will be denoted by
$d_k.$ Two chains of cuts $\{\sigma_m\}$ and
$\{\sigma_k^{\,\prime}\}$ are called {\it equivalent} if for each $m
=1,2, \ldots $ the domain $d_m$ contains all domains $d_k^{\,\prime}
$ except for a finite number, and for each $k= 1,2, \ldots $, the
domain $d_k^{\,\prime} $ also contains all domains $d_m$ for
excluding a finite number.

\medskip
The following statement contains some information on important
properties of domains with locally quasiconformal boundaries.

\begin{lemma}\label{pr7}
{\sl\, Suppose that a domain $D$ in ${\Bbb S} $ has a locally
quasiconformal boundary, $p_0\in \partial D$ and $r_m> 0$ is an
arbitrary sequence such that $r_m \rightarrow 0$ as $ m\rightarrow
0. $ Let $ U $ be a neighborhood of the point $ p_0, $ for which
there is a quasiconformal mapping $\psi: U \rightarrow {\Bbb D},$
$\psi(\partial D \cap U)={\Bbb D}\cap I,$ $I:=(- 1, 1)=\{z \in {\Bbb
D}: x=0, z=x+iy \}.$ Then:

\medskip
1)  there is a sequence of neighborhoods $U_k,$ $k=1,2, \ldots, $ of
the point $p_0,$ contracting to this point, such that
$\psi(U_k)=B(0, 1/2^k),$ $\psi(U_k\cap D)=\{z=x+iy\in B(0, 1/2^k):
x>0\};$

\medskip
2) the sets $\sigma_k:=\partial U_k\cap D$ are cuts of the domain
$D,$ in this case, $\widetilde{h}(\sigma_k)\rightarrow 0$ as
$k\rightarrow\infty,$ $\widetilde{h}(\sigma_k):=\sup\limits_{x,
y\in\sigma_k}\widetilde{h}(x, y),$ and the corresponding domain
$d_k$ is $U_k\cap D;$

\medskip
3) there is a subsequence $r_ {m_l}\rightarrow 0,$ $l\rightarrow
\infty,$ and the corresponding sequence of arcs $\gamma_l\subset
\widetilde{S}(p_0, r_ {m_l})\cap D,$ $l= 1,2, \ldots, $ forming a
chain of cuts equivalent to $\sigma_k, $ $k= 1,2, \ldots \,. $ }
\end{lemma}

\begin{proof}
Arguing similarly to the proof of Theorem~17.10 in~\cite{Va}, we
show that the neighborhood $U$ in the definition of a locally
quasiconformal boundary can be chosen arbitrarily small. Indeed, by
definition, the point $p_0\in \partial D $ has a neighborhood $U,$
which can be mapped by a quasiconformal mapping $\psi$ onto the unit
disk ${\Bbb D}\subset{\Bbb C}$ so that $\psi(\partial D \cap U)=I,$
where $I:=(- 1, 1)=\{z\in {\Bbb D}: x=0, z=x+iy\}.$ Since $\psi$ is
a homeomorphism, then either $\psi(U\cap D)={\Bbb D}_+,$ or
$\psi(U\cap D)={\Bbb D}\setminus \overline{{\Bbb D}_ +},$ where
${\Bbb D}_+:=\{z\in {\Bbb D}: x> 0, z=x+iy\}.$ Thus, without loss of
generality, we may assume that $\psi(U \cap D)={\Bbb D}_+.$ Choose a
neighborhood $V\subset U $ containing the point $p_0.$ If
$r<1-|\psi(p_0)|,$ then, by the triangle inequality, the ball
$B(\psi(p_0),r)$ lies strictly inside ${\Bbb D}.$ Since $\psi$ is a
homeomorphism in $U,$ then, in particular, $\psi^{\,- 1}$ is a
continuous mapping. In this case, there is $r<1-|\psi(p_0)| $ with
the following property: the condition $|\psi(p)-\psi(p_0)|<r$
implies that $p\in V.$ In addition, if $\psi(p)\in{\Bbb D}_+\cup I$
and $|\psi(p)-\psi(p_0)|<r,$ then $p\in V \cap \overline{D}.$
Setting $U_1:=\psi^{\,- 1}(B(\psi(p_0), r)),$ we note that
$U_1\subset V$ and $U_1$ is a neighborhood of the point $p_0.$ In
this case, $U_1\cap \overline{D}=\psi^{\,-
1}(B(\psi(p_0),r)\cap({\Bbb D}_+\cup I)).$ Setting
$H(p)=(\psi(p)-\psi(p_0))/r,$ we obtain the mapping $H$ of the
neighborhood $U_1$ on ${\Bbb D}$ such that
$H(U_1\cap\overline{D})={\Bbb D}_+\cup I$ and $H(p_0)=0.$ Since $H$
is a homeomorphism, it follows that $H(U_1\cap \partial D)=I.$ It is
also clear that if the original mapping $\psi$ is quasiconformal,
then the same is the mapping $H.$ Thus, the neighborhood $U_1$
satisfies all the same conditions as the original neighborhood $U.$
In what follows, we use the notation $\psi $ instead of $H,$ and we
assume that $\psi(p_0)=0.$

\medskip
From the above reasoning it follows that there is a decreasing
sequence of neighborhoods $U_k$ of the point $p_0,$ for which
$p_0=\bigcap\limits_{k=1}^{\infty}\overline{U_k\cap D},$
$\psi(U_k)=B(0, 1/2^k),$ $\psi(\partial U_k\cap D)=S(0, 1/2^k)\cap
{\Bbb D}_+.$ By direct calculations it is easy to see that that the
sequence $\sigma_k:=\partial U_k \cap D$ forms a chain of cuts of
the domain $D.$ From the equality
$\bigcap\limits_{k=1}^{\infty}\overline{U_k\cap D}=p_0$  it follows
that $\widetilde{h}(\sigma_k)\leqslant
\widetilde{h}(\overline{U_k\cap D})\rightarrow 0$ as
$k\rightarrow\infty,$ where we use the notation
$\widetilde{h}(A):=\sup\limits_{x, y\in A}\widetilde{h}(x, y).$
Items~1) and~2) of Lemma~\ref{pr7} are established. It remains to
establish item~3).
\begin{figure}[h]
\centerline{\includegraphics[scale=0.52]{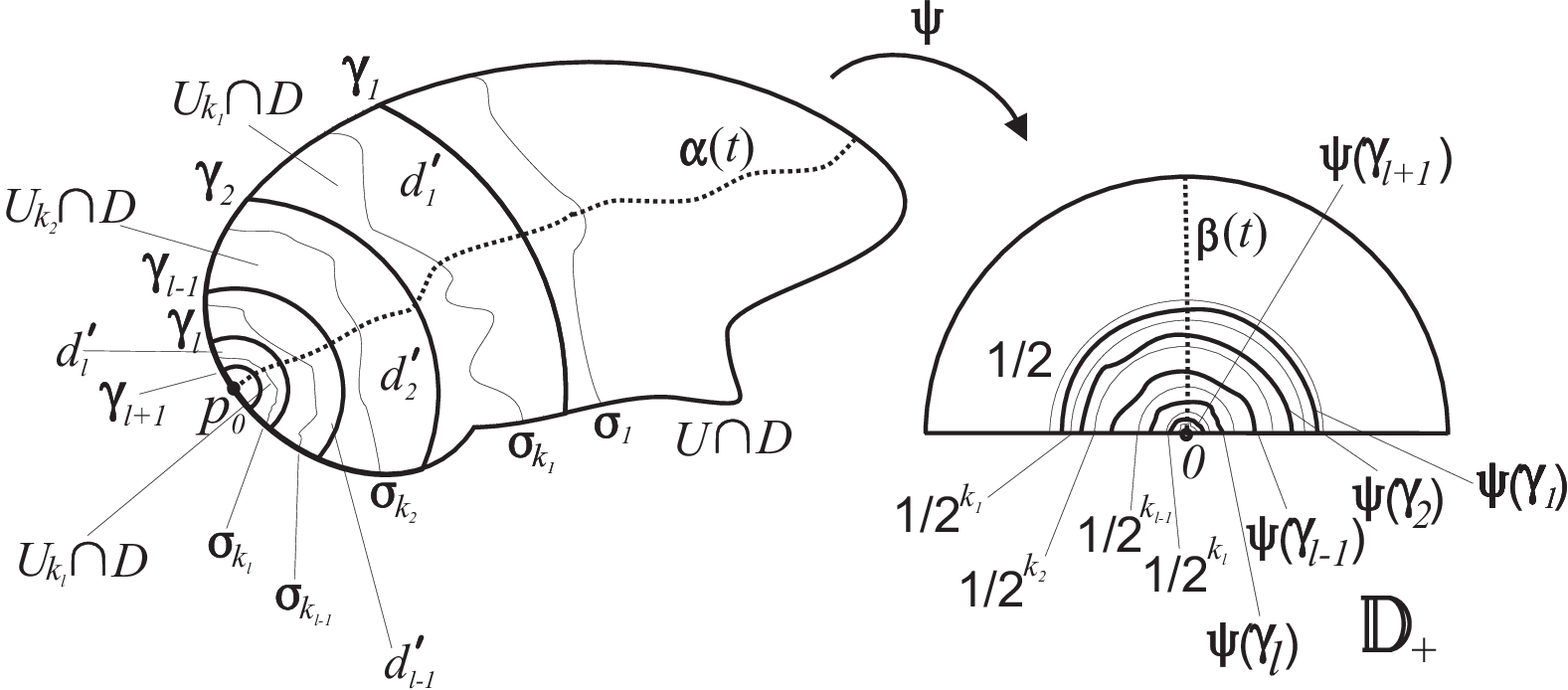}} \caption{On
the proof of Lemma~\ref{pr7}}\label{fig3}
\end{figure}
Consider the segment $\beta(t)=it,$ $t\in (0, 1),$ in ${\Bbb D}_+.$
Put $\alpha(t):=\psi^{\,-1}(\beta(t)).$ Then $\alpha$ is a path in
$U\cap D$ with origin at the point $p_0.$ Let $m_1\in {\Bbb N} $ be
that $r_{m_1}<\widetilde{h}(p_0, \sigma_1).$ By~\cite[Theorem~1.I.5,
$\S \,46$]{Ku} $\widetilde{S}(p_0, r_{m_1})\cap
\alpha\ne\varnothing.$ Let $\gamma_1$ be an arc of the set
$\widetilde{S}(p_0, r_{m_1})\cap D,$ such that $\alpha(t_1)\in
\gamma_1,$ where $t_1:=\min\{t\in (0, 1): \alpha(t)\in
\widetilde{S}(p_0, r_{m_1})\}$ (see~Figure~\ref{fig3}).

By construction, the ends of the path $\gamma_1$ lie on $\partial
D,$ therefore $\psi(\gamma_1)$ is a path whose ends lie on the
segment $I\subset{\Bbb D}.$ Obviously, $\psi(\gamma_1)$ splits
${\Bbb D}_+$ into two domains. Therefore, $\gamma_1$ divides the
domain $D$ into two domains, as well. Let $d_1^{\,\prime}$ be the
component of $D\setminus \gamma_1,$ containing the path
$\alpha_1:=\alpha|_{(0, t_1)}.$ Note that there is $k_1\in {\Bbb N}$
such that $\overline{U_ {k_1}}\cap D \subset d_1^{\,\prime}.$
Indeed, by construction $\overline{U_{k_1}}\subset
\widetilde{B}(p_0, r_{m_1})$ for sufficiently large $k_1\in{\Bbb N}$
and the point $\alpha_1(t)$ belongs to $U_{k_1}\cap D $ for
sufficiently small $t,$ since $U_{k_1}$ is a neighborhood of the
point $p_0.$ Hence, $\overline{U_ {k_1}}\cap D $ belongs to some
component of $D\setminus \gamma_1$ containing $\alpha_1,$ that is,
$\overline{U_{k_1}}\cap D \subset d_1^{\,\prime}.$ Note that
$\sigma_{k_1}\subset d_1^{\,\prime},$ because $\sigma_{k_1}=\partial
U_{k_1}\cap D \subset d_1^{\,\prime}.$

\medskip
Observe also that $d_1^{\,\prime}\subset D\setminus \sigma_1.$ This
follows from the fact that $d_1^{\,\prime}$ is bounded by $\gamma_1$
and some part of the boundary of the domain $U\cap D,$ and
$\sigma_1$ lies in the domain bounded by $\gamma_1$ and another part
of the boundary $U\cap D.$ Thus, $\gamma_1$ separates $\sigma_1$
from $U_{k_1}\cap D$ in $D; $ in particular, $\gamma_1$ separates
$\sigma_1$ from $\sigma_{k_1}$ in $D.$

\medskip
Let $m_2 \in{\Bbb N},$ $m_2>m_1,$ be such that
$r_{m_2}<\widetilde{h}(p_0, \sigma_{k_1}).$ By~\cite[Theorem~1.I.5,
$\S\,46$]{Ku} $\widetilde{S}(p_0, r_{m_2})\cap
\alpha\ne\varnothing.$ Let $\gamma_2$ be an arc of the set
$\widetilde{S}(p_0, r_{m_2})\cap D$ such that $\alpha(t_2)\in
\gamma_2,$ where $t_2:=\min\{t\in (0, 1): \alpha(t)\in
\widetilde{S}(p_0, r_{m_2})\}.$ By construction, the ends of the
path $\gamma_2$ lie on $\partial D,$ therefore $\psi(\gamma_2)$ is a
path whose ends lie on the segment $I\subset{\Bbb D}.$ Obviously,
$\psi(\gamma_2)$ splits ${\Bbb D}_+$ into two domains, therefore,
$\gamma_2$ divides the domain $D$ into two domains. Let
$d_2^{\,\prime}$ be a component of $D\setminus \gamma_2,$ containing
a path $\alpha_2:=\alpha|_{(0, t_2)}.$ Observe that there is
$k_2\in{\Bbb N}$ such that $\overline{U_{k_2}}\cap D\subset
d_2^{\,\prime}.$ Indeed, by construction $\overline{U_{k_2}}\subset
\widetilde{B}(p_0, r_{m_2})$ for large enough $k_2\in {\Bbb N}$ and
$\alpha_2(t)$ belongs to $U_{k_2}\cap D$ for sufficiently small $t,$
since $U_{k_2}$ is a neighborhood of the point $p_0.$ Hence,
$\overline{U_{k_2}}\cap D$ belongs to some component of $D\setminus
\gamma_2,$ containing $\alpha_2,$ that is, $\overline{U_{k_2}}\cap
D\subset d_2^{\,\prime}.$ Observe that $\sigma_{k_2}=\partial
U_{k_2}\cap D\subset d_2^{\,\prime}.$

\medskip
Observe also that, $d_2^{\,\prime}\subset D\setminus \sigma_{k_1}.$
Thus $\gamma_2$ separates $\sigma_{k_1}$ from $U_{k_2}\cap D$ in
$D;$ in particular, $\gamma_2$ separates $\sigma_{k_1}$ from
$\sigma_{k_2}$ in $D.$

\medskip
For the same reason, $\gamma_1\subset D\setminus d_2^{\,\prime}.$
Indeed, $d_2^{\,\prime}$ is one of the components of $D\setminus
\gamma_1,$ not the same as $d_1^{\,\prime},$ in particular,
$d_2^{\,\prime}\subset D\setminus \gamma_1.$ Thus, $\gamma_2 $ also
separates $\gamma_1$ from $U_{k_2}\cap D$ in $D;$ in particular,
$\gamma_2$ separates $\gamma_1$ from $\sigma_{k_2}$ in $D.$

\medskip
Etc. As a result of the endless process, we will have that some
sequence of arcs $\gamma_l\subset \widetilde{S}(p_0, r_{m_l}),$
$l=1,2, \ldots,$  separating $\sigma_{k_l}$ from $\sigma_{k_{l-1}}$
in $D,$ in this case,

\medskip
\textbf{(1)} $U_{k_l}\cap D\subset d_l^{\,\prime}$
and, in addition,

\medskip
\textbf{(2)} $\gamma_l$ separates $\gamma_{l-1}$ from $\sigma_{k_l}$
in $D,$

\medskip
\textbf{(3)} $\sigma_{k_l}\in d_l^{\,\prime};$

\medskip
\textbf{(4)} $\gamma_{l+1}\subset U_{k_l}\cap D.$

\medskip
Let us show that $\gamma_l$ separates $\gamma_{l+1}$ from
$\gamma_{l-1}$ for any $l\in{\Bbb N},$ more precisely, show that
$\gamma_{l+1}\subset d^{\,\prime}_l$ and $\gamma_{l-1}\subset
D\setminus d^{\,\prime}_l.$ Indeed, step by step~\textbf{(1)}
and~\textbf{(4)} $\gamma_{l+1}\subset U_{k_l}\cap D\subset
d_l^{\,\prime}.$ Since as proved, $\gamma_l$ separates
$\gamma_{l-1}$ from $\sigma_{k_l}$ in $D$ and $\sigma_{k_l}\subset
d_l^{\,\prime}$ by~\textbf{(3)}, by~\cite[Theorem~1.I.5,
$\S\,46$]{Ku} $\gamma_{l-1}\subset D\setminus d_l^{\,\prime},$ as
required to establish.

\medskip
It follows from the above that the sequence of cuts $\gamma_l,$
$l=1,2,\ldots ,$ forms a chain. Note that the sequences of cuts
$\gamma_l$ and $\sigma_k$ are equivalent. Indeed, given $l\in{\Bbb
N},$ by the property~\textbf{(1)} $U_{k_l}\cap D\subset
d_l^{\,\prime},$ therefore also $U_k\cap D\subset d_l^{\,\prime}$
for $k\geqslant k_l.$ Conversely, fix $k \in{\Bbb N}$ and consider
the corresponding number $l=l(k)\in {\Bbb N}$ such that
$k_l\geqslant k,$ where $k_l,$ $l=1,2,\ldots $ is the subsequence
constructed above. Notice, that $d_{l+1}^{\prime}$ belongs to
exactly one of the components of $D\setminus \sigma_{k_l},$ namely,
either $d_{l+1}^{\,\prime}\subset U_{k_l}\cap D,$ or
$d_{l+1}^{\,\prime}\subset D\setminus U_{k_l}.$ On the other hand,
by the condition~\textbf{(3)} $\sigma_{k_{l+1}}\subset
d_{l+1}^{\,\prime},$ in addition, $\sigma_{k_{l+1}}\subset
U_{k_{l+1}-1}\cap D\subset U_{k_l}\cap D,$ because $\sigma_k,$
$k=1,2,\ldots $ is a chain of cuts. Moreover, $k_{l+1}-1\geqslant
k_l.$ In this case, $d_{l+1}^{\,\prime}$ belongs to the component of
$D\setminus\sigma_{k_l},$ containing $\sigma_{k_{l+1}},$ that is,
$d_l^{\,\prime}\subset U_{k_l}\cap D\subset U_k\cap D.$ Then also
$d_l^{\,\prime}\subset U_k\cap D$ for any $l\geqslant l(k).$
Equivalence of chains $\gamma_l$ and $\sigma_k$ is established.
Lemma~\ref{pr7} is completely proved.~$\Box$
\end{proof}

\medskip
The following lemma is technically necessary to establish the main
result on the boundary behavior of mappings. We have specially
highlighted it into a separate statement, emphasizing that it refers
to mappings satisfying rather general topological conditions.

\begin{lemma}\label{lem6}
{\sl\, Let $D$ and $D_{\, *}$ be domains on Riemannian surfaces
${\Bbb S}$ and ${\Bbb S}_*,$ respectively, let $\overline{D_{\,*}}$
be a compactum in ${\Bbb S}_*,$ let $p_0\in\partial D$ and let
$Q:D\rightarrow (0, \infty)$ be a given function measurable with
respect to measure $\widetilde{h}.$ Let also $f:D\rightarrow
D_{\,*}$ be an open discrete closed mapping of the domain $D$ onto
the domain $f(D)=D_{\,*}.$ Suppose the domain $D$ has a locally
quasiconformal boundary , and the boundary of the domain $D_{\,*}$
is strongly accessible.

Suppose that $p_0\in \partial D, $ and that there are at least two
sequences $p_i,$ $p_i^{\,\prime}\in D,$ $i=1,2,\ldots,$ such that
$p_i\rightarrow p_0,$ $p_i^{\,\prime}\rightarrow p_0$ as
$i\rightarrow \infty,$ $f(p_i)\rightarrow y,$
$f(p_i^{\,\prime})\rightarrow y^{\,\prime}$ as $i\rightarrow \infty$
and $y^{\,\prime}\ne y.$

Then there are  $0<\delta_0^{\,\prime}$ and $l_0>0$  such that the
inequality
\begin{equation}\label{eq42}
l(f(\widetilde{S}(p_0, r)\cap D))\geqslant l_0,\quad\forall\,\,r\in
(0, \delta_0^{\,\prime})\,,
\end{equation}
where $l$ denotes the length of the path (dashed line) on the
Riemannian surface ${\Bbb S}_*. $ }
\end{lemma}

\begin{proof}
By the definition of a strongly accessible boundary at the point
$y\in \partial D_{\,*},$ for the neighborhood $U$ of the point $y,$
not containing the point $y^{\,\prime},$ there is a compact set
$C_0^{\,\prime}\subset D_{\,*},$ a neighborhood $V$ of the point
$y,$ $V\subset U,$ and a number $\delta> 0$ such that
\begin{equation}\label{eq1AB}
M(\Gamma(C_0^{\,\prime}, F, D_{\,*}))\ge \delta
>0
\end{equation}
for an any continuum $F,$ intersecting $\partial U $ and $\partial
V. $ By Lemma~\ref{pr7}, there is a sequence of neighborhoods $U_i,$
$i=1,2, \ldots, $ of $p_0,$ such that the set $d_i:=U_i\cap D$ is
connected. Without loss of generality, we may assume that $p_i$ and
$p_i^{\,\prime}$ belong to $d_i.$ In this case, join the points
$p_i$ and $p_i^{\,\prime}$ by the path $\alpha_i,$ lying in $d_i.$
Since $f(p_i)\in V $ and $f(p_i^{\,\prime})\in D\setminus
\overline{U}$ for sufficiently large $i\in {\Bbb N},$ there is a
number $i_0 \in {\Bbb N},$ such that by~(\ref{eq1AB})

\begin{equation}\label{eq2D}
M(\Gamma(C_0^{\,\prime}, f(|\alpha_i|), D_{\,*}))\ge \delta
>0
\end{equation}
for any $i\ge i_0\in {\Bbb N}$ (see Figure~\ref{fig2}).
\begin{figure}[h]
\centerline{\includegraphics[scale=0.5]{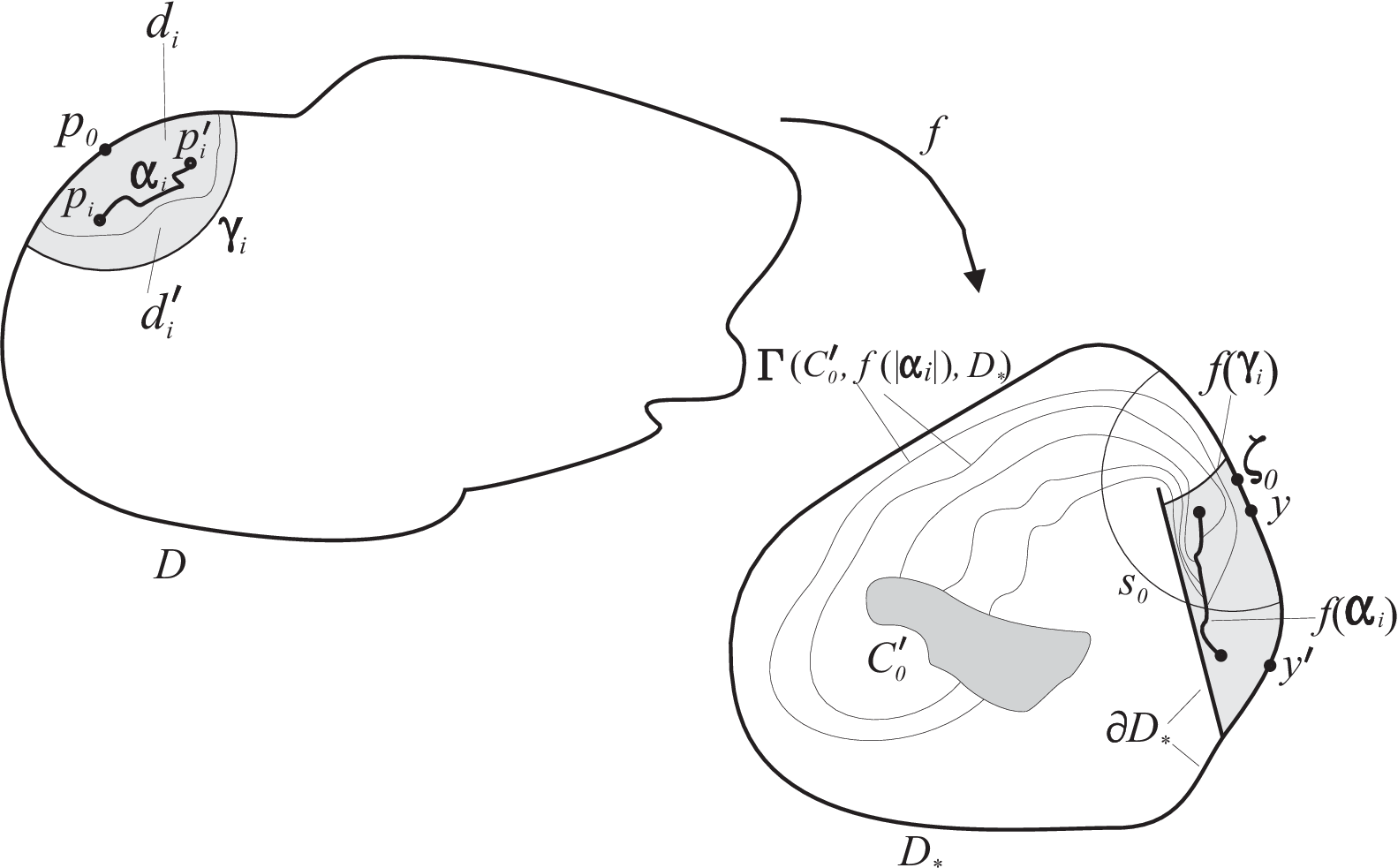}} \caption{To
the proof of Lemma~\ref{lem6}}\label{fig2}
\end{figure}
Let us prove Lemma~\ref{lem6} by contradiction. Suppose
that~(\ref{eq42}) does not hold. Then for any $k\in {\Bbb N}$ there
is $r=r_k>0$ such that $l(f(\widetilde{S}(p_0, r_k)\cap D))< 1/k,$
$r_k\rightarrow 0,$ $k\rightarrow\infty.$
By Lemma~\ref{pr7}, there is a subsequence $r_{k_i}$ of the sequence
$r_k$ and some sequence of arcs $\gamma_i\subset \widetilde{S}(p_0,
r_ {k_i})$ such that $d_i\subset d^{\,\prime}_i$ and
$d^{\,\prime}_i$ is the corresponding component of the set
$D\setminus |\gamma_i|,$ $i=1,2,\ldots .$ Let $\zeta_i,$
$i=1,2,\ldots ,$ is an arbitrary sequence of points from
$f(|\gamma_i|).$ Since $\overline{D_{\, *}}$ is a compactum in
${\Bbb S}_*,$ we may assume that $\zeta_i\rightarrow \zeta_0$ as
$i\rightarrow\infty,$ $\zeta_0\in \overline{D_{\,*}}.$
Observe that $\zeta_i=f(x_i),$ $x_i\in \widetilde{S}(p_0,
r_{k_i})\cap D,$ so $\zeta_0\in
\partial D_{\,*}$ by the closeness of the mapping $f$ and
Proposition~\ref{pr4}.

\medskip
Note that, since the mapping $f$ is closed, there is a number
$i_0\in {\Bbb N}$ such that
\begin{equation}\label{eq2C}
C_0^{\,\prime}\subset D_*\setminus \overline{f(d^{\,\prime}_i)}\,.
\end{equation}
for $i\geqslant i_0.$
Indeed, if we assume that the inclusion~(\ref{eq2C}) fails, then
there is a sequence $i_m>0,$ $m=1,2,\ldots,$ $i_m\rightarrow \infty$
as $m\rightarrow\infty,$ and $y_m\in
\overline{f(d^{\,\prime}_{i_m})}\cap C_0^{\,\prime}.$ Since
$C_0^{\,\prime}$ is a compactum in $f(D),$ we may assume that
$y_m\rightarrow y_0\in C_0^{\,\prime} $ as $m\rightarrow\infty.$
Since~$y_m\in \overline{f(d^{\,\prime}_{i_m})}\cap C_0^{\,\prime},$
for any $m\in{\Bbb N}$ there is a sequence $y_{km}\in
\overline{f(d^{\,\prime}_{i_m})}$ such that $y_{km}\rightarrow y_m$
as $k\rightarrow\infty.$ Observe that $y_{km}=f(q_{km}),$ $q_{km}\in
d^{\,\prime}_{i_m}.$

Due to the convergence of $y_{k1}$ to $y_1,$ for the number $1/2,$
there is a number $k_1$ such that ${\widetilde{h_*}}(y_1,
y_{k_11})<1/2.$ Similarly, due to the convergence of $y_{k1}$ to
$y_2,$ for the number $1/4$ there is a number $k_2$ such that
${\widetilde{h_*}}(y_2, y_{k_22})<1/4.$ Generally, since $y_{km}$
converges to $y_m,$ for the number $1/2^ m$ there is a number $k_m$
such that ${\widetilde{h_*}}(y_m, y_{k_mm})<1/2^m.$ But then, since
by the construction of $y_m \rightarrow y_0$ as $m\rightarrow
\infty, $ for any fixed $\varepsilon>0$ we have that
$${\widetilde{h_*}}(y_0,y_{k_mm})\leqslant {\widetilde{h_*}}(y_0, y_m)
+{\widetilde{h_*}}(y_m, y_{k_mm})\leqslant \varepsilon+1/2^m$$
for any $m\geqslant M=M(\varepsilon),$ and thus $y_{k_mm}\rightarrow
y_0$ as $m\rightarrow\infty.$ But, on the other hand hand,
$y_{k_mm}=f(q_{k_mm}),$ $q_{k_mm}\in d_{i_m}^{\,\prime},$ therefore
$y_0\in C(f, p_0),$ which contradicts the closeness of the mapping
$f.$ Indeed, by Lemma~\ref{pr7} there exists a sequence $t_i> 0,$
$t_i\rightarrow 0$ as $i\rightarrow\infty,$ such that
$d_i^{\,\prime}\subset \widetilde{B}(p_0, t_i).$ Then $y_0 \in C (f,
\partial D)\subset \partial D_*$ (see Proposition~\ref{pr4}).
At the same time, $y_0\in C_0^{\,\prime},$ that is, $y_0$ is an
inner point of the domain $D_*.$ The resulting contradiction
indicates the validity of the inclusion~(\ref{eq2C}).

\medskip
Let us now show that
\begin{equation}\label{eq7A}
\partial f(d^{\,\prime}_i)\cap f(D)\subset
f(|\gamma_i|)
\end{equation}
for any $i\in{\Bbb N}.$

Indeed, let $y_0\in \partial f(d^{\,\prime}_i)\cap f(D),$ then there
is a sequence $y_m\in f(d^{\,\prime}_i)\cap f(D)$ such that $y_m\in
f(d^{\,\prime}_i)\cap f(D),$ $y_m\rightarrow y_0$ as $m\rightarrow
\infty,$ where $y_m=f(\xi_m),$ $\xi_m\in d^{\,\prime}_i.$ Without
loss of generality, we may assume that $\xi_m\rightarrow \xi_0$ as
$m\rightarrow\infty.$ Note that the case $\xi_0\in\partial D $ is
impossible, since in this case $y_0\in C(f,
\partial D),$ which contradicts the closeness of the mapping $f.$
Then $\xi_0 \in D.$ Two situations are possible: 1) $ \xi_0 \in
d^{\,\prime}_i$ and 2) $\xi_0 \in |\gamma_i |.$ Note that case 1) is
impossible, since, in this case, $f(\xi_0)= y_0$ and $y_0$ is an
inner point of the set $f(d^{\,\prime}_i)$ by the openness of the
mapping $f,$ which contradicts the choice of $y_0.$ Thus, the
inclusion~(\ref{eq7A}) is established.

\medskip
By the assumption $l(f(\widetilde{S}(p_0, r_k)\cap D))< 1/k,$
$r_k\rightarrow 0$ as $k\rightarrow\infty,$ we obtain that
$$\widetilde{h}_*(f(|\gamma_i|))\rightarrow 0$$ as
$i\rightarrow\infty,$
$\widetilde{h}_*(f(|\gamma_i|)):=\sup\limits_{p_*, q_*\in
f(|\gamma_i|)}\widetilde{h}_*(p_*, q_*).$
Then, for any $s\in{\Bbb N} $ there is a number $i_s\in{\Bbb N}$
such that $f(|\gamma_{i_s}|)\subset \widetilde{B}(\zeta_0, 1/s).$
Since $C_0^{\,\prime}$ is a compactum in $D_{\,*},$ there is
$s_0\geqslant 1$ such that $C_0^{\,\prime}\cap
\widetilde{B}(\zeta_0, 1/{s_0})=\varnothing.$

\medskip
Now, note that $\Gamma(f(|\gamma_{i_s}|), C_0^{\,\prime},
D_{\,*})>\Gamma(\widetilde{S}(\zeta_0, 1/{s_0}),
\widetilde{S}(\zeta_0, 1/s), D_{\,*})$ for any $s>s_0$
(see~\cite[Theorem~1.I.5, $\S\,46$]{Ku}). Hence, by the minorization
of the modulus of families of paths and by Proposition~\ref{pr5}
\begin{equation}\label{eq43}
M(\Gamma(f(|\gamma_{i_s}|), C_0^{\,\prime}, D_{\,*}))\leqslant
M(\Gamma(\widetilde{S}(\zeta_0, 1/{s_0}), \widetilde{S}(\zeta_0,
1/s), D_{\,*}))\rightarrow 0
\end{equation}
as $s\rightarrow\infty.$ On the other hand, recall that
$|\alpha_{i_s}|\subset d_{i_s}\subset d^{\,\prime}_{i_s}.$ Now,
by~(\ref{eq2C})
$$|\beta|\cap f(d^{\,\prime}_{i_s})\ne\varnothing\ne |\beta|\cap (D_{\,*}\setminus
f(d^{\,\prime}_{i_s})\,,$$
for any path $\beta\in\Gamma(f(|\alpha_{i_s}|), C^{\,\prime}_0,
D_{\,*}).$ Thus, by~\cite[Theorem~1.I.5, $\S\,46$]{Ku} and
by~(\ref{eq7A}),
$$\Gamma(f(|\alpha_{i_s}|), C^{\,\prime}_0, D_{\,*})> \Gamma(f(|\gamma_{i_s}|), C^{\,\prime}_0,
D_{\,*}).$$
From this and by~(\ref{eq43}), we obtain
$M(\Gamma(f(|\alpha_{i_s}|), C^{\,\prime}_0, D_{\,*}))\rightarrow 0$
as $s\rightarrow \infty,$ which contradicts relation~(\ref{eq2D}).
The resulting contradiction indicates the validity of the
inequality~(\ref{eq42}).~$\Box$
\end{proof}

\begin{theorem}\label{th1}{\sl\,
Let $D$ and $D_{\,*}$ be domains on Riemannian surfaces ${\Bbb S}$
and ${\Bbb S}_*,$ respectively, let $\overline{D_{\,*}}$ be a
compactum in ${\Bbb S}_*,$ let $p_0\in\partial D $ and let
$Q:D\rightarrow (0, \infty)$ be a given function measurable with
respect to the measure $\widetilde{h}.$ Let also $f:D\rightarrow
D_{\,*}$ be an open discrete closed mapping of $D$ onto
$f(D)=D_{\,*}$ with the condition~(\ref{eq1A}) at the point $p_0.$
Suppose that the domain $D$ has a locally quasiconformal boundary,
and the boundary of the domain $D_{\,*}$ is strongly accessible. If
the relations
\begin{equation}\label{eq9AA}
\int\limits_{\varepsilon}^{\varepsilon_0} \frac{dt}{\Vert
Q\Vert(t)}<\infty\,,\qquad \int\limits_0^{\varepsilon_0}
\frac{dt}{\Vert Q\Vert(t)}=\infty\,,
\end{equation}
hold for some $0<\varepsilon_0$ and any $\varepsilon\in (0,
\varepsilon_0),$ where $\Vert
Q\Vert(r):=\int\limits_{\widetilde{S}(p_0,
r)}Q(p)\,ds_{\widetilde{h}}(p),$ then $f$ has a continuous extension
to~$p_0.$}
\end{theorem}

\begin{proof}
Suppose the opposite. Then, since $\overline{D_{\,*}}$  is compact
in ${\Bbb S}_*,$ there are at least two sequences $p_i,$
$p_i^{\,\prime}\in D,$ $i=1,2,\ldots,$ such that $p_i\rightarrow
p_0,$ $p_i^{\,\prime}\rightarrow p_0$ as $i\rightarrow \infty,$
$f(p_i)\rightarrow y,$ $f(p_i^{\,\prime})\rightarrow y^{\,\prime}$
as $i\rightarrow \infty$ and $y^{\,\prime}\ne y.$ Let
$\Gamma^{\,\delta}_i$ be a family of all dished lines
$f(\widetilde{S}(p_0, r)\cap D),$ $r\in (2^{\,-i}, \delta).$ By
Lemmas~\ref{lem4A} and~\ref{lem2}, there is $0<d_0<\sup\limits_{p\in
D}\widetilde{h}(p, p_0)$ such that
\begin{equation}\label{eq8A}
M(\Gamma^{\,\delta_0}_i)\geqslant \int\limits_{2^{\,-i}}^{\delta_0}
\frac{dr}{\Vert\,Q\Vert(r)}\quad\forall\,\, i\in {\Bbb
N}\,,\end{equation} for any $0<\delta_0<d_0,$ where
$\Vert
Q\Vert(r)=\int\limits_{D(p_0,r)}Q(p)\,\,ds_{\widetilde{h}}(p)$
denotes $L_1$-norm of the function $Q$ under the circle
$D(p_0,r):=\widetilde{S}(p_0,r)\cap D.$ By~(\ref{eq8A})
and~(\ref{eq9AA}) we obtain that
\begin{equation}\label{eq40A}
M(\Gamma^{\,\delta_0}_i)\rightarrow\infty\,,\quad
i\rightarrow\infty\,.
\end{equation}
On the other hand, by Lemma~\ref{lem6} there are
$0<\delta_0^{\,\prime}<d_0$ and $l_0>0$ such that
$$
l(f(\widetilde{S}(p_0, r)\cap D))\geqslant l_0,\quad\forall\,\,r\in
(0, \delta_0^{\,\prime})\,,
$$
where $l$ denotes the length of the dished line in ${\Bbb S}_*.$ In
particular, the function
$$\rho(p)= \left \{\begin{array}{rr} 1/l_0\ , & \ p\in D_{\,*}\ ,
\\ 0\ ,  &  p\not\in D_{\,*}
\end{array} \right.$$
is admissible for $\Gamma_i^{\delta^{\,\prime}_0},$
$0<\delta_0^{\,\prime}<d_0.$ Since $\overline{D_{\,*}}$ is a
compactum in ${\Bbb S}_*,$ the $\widetilde{h_*}$-area of $D_{\,*}$
is finite. Therefore,
$$M(\Gamma_i^{\delta^{\,\prime}_0})\leqslant 1/l_0^2\cdot \widetilde{h_*}(D_{\,*})<\infty\,.$$
The last condition contradicts the relation~(\ref{eq40A}) for
$\delta_0:=\delta^{\,\prime}_0,$ which refutes the assumption that
the mapping $f$ has no limit at the point~$p_0.$~$\Box$
\end{proof}

{\it Proof of Theorem~\ref{th2}} follows immediately from
Theorems~\ref{thOS4.2} and \ref{th1}.~$\Box$

\medskip Let $p_0\in {\Bbb S}$ and let $\varphi:{\Bbb
S}\rightarrow {\Bbb R}$ be a function integrable in some
neighborhood $U$ of the point $p_0$ with respect to $\widetilde{h}.$
Following~\cite[Section~2]{IR} (see also~\cite[Section~6.1,
Ch.~6]{MRSY}), we say that a function $\varphi:{\Bbb S}\rightarrow
{\Bbb R}$ has a {\it finite mean oscillation} at the point $p_0\in
D$, we write $\varphi\in FMO (p_0),$ if
%
%
%
%\begin{equation}\label{eq29*!}
%
$$\limsup\limits_{\varepsilon\rightarrow
0}\frac{1}{\widetilde{h}(\widetilde{B}(p_0,
\varepsilon))}\int\limits_{\widetilde{B}(p_0,\,\varepsilon)}
|{\varphi}(p)-\overline{\varphi}_{\varepsilon}|\
d\widetilde{h}(p)<\infty\,,$$
%
%\end{equation}
%
where
$\overline{{\varphi}}_{\varepsilon}=\frac{1}
{\widetilde{h}(\widetilde{B}(p_0,
\varepsilon))}\int\limits_{\widetilde{B}(p_0, \varepsilon)}
{\varphi}(p) \,d\widetilde{h}(p).$ In what follows, we will talk
about results related to the function of the finite mean
oscillation, therefore it is extremely important for us to use the
following two most important facts related to these functions.

Let $D$ be a domain in ${\Bbb S}, $ and let $\varphi: {\Bbb S}
\rightarrow {\Bbb R}$ be a nonnegative function with a finite mean
oscillation at the point $p_0\in \overline{D}\subset {\Bbb S},$
$\varphi (x)=0 $ for $x\not \in D.$ By~\cite[Theorem~7.2.2]{Berd},
the surface ${\Bbb S}$ is locally Ahlfors 2-regular, so that
by~\cite[Lemma~3]{Sm}
\begin{equation}\label{eq44}
\int\limits_{\varepsilon<\widetilde{h}(p, p_0)<
\varepsilon_0}\frac{\varphi(p)\, d\widetilde{h}(p)}
{\left(\widetilde{h}(p, p_0)\log\frac{1}{\widetilde{h}(p,
p_0)}\right)^2} = O \left(\log\log \frac{1}{\varepsilon}\right)
\end{equation}
as $\varepsilon \rightarrow 0$ for some $0<\varepsilon_0<{\rm
dist}(p_0, \partial U) $ and some normal neighborhood $U$ of the
point $p_0.$ The following statement may be proved similarly
to~\cite[Lemma~7.4, Ch.~7]{MRSY}, cf.~\cite[Lemma~3.7]{RS$_2$}
or~\cite[Lemma~4.2]{ARS}.

\begin{proposition}\label{pr1A}
{\sl\, Let $p_0 \in {\Bbb S},$ let $U$ be some normal neighborhood
of $p_0,$ $0<r_1<r_2<{\rm dist}\,(p_0,
\partial U),$ and let $Q:{\Bbb S}\rightarrow [0,
\infty]$ be an integrable function in $U$ with respect to the
measure $\widetilde{h}.$ Set $\widetilde{A}=\widetilde{A}(p_0, r_1,
r_2)=\{p\in {\Bbb S}: r_1<\widetilde{h}(p, p_0)<r_2\},$ $\Vert
Q\Vert(r)=\int\limits_{\widetilde{S}(p_0,
r)}Q(p)\,\,ds_{\widetilde{h}}(p),$
$\eta_0(r):=\frac{1}{J\cdot\Vert Q\Vert(r)},$
where $J=J(p_0,r_1,r_2):=\int\limits_{r_1}^{r_2}\ \frac{dr}{\Vert
Q\Vert(r)}\,.$
Then
$$J^{\,-1}=\int\limits_{\widetilde{A}(p_0, r_1, r_2)} Q(p)\cdot
\eta_0^2(\widetilde{h}(p, p_0))\ d\widetilde{h}(p)\leqslant$$
\begin{equation}\label{eq10A}
\leqslant\int\limits_{\widetilde{A}(p_0, r_1, r_2)} Q(p)\cdot
\eta^2(\widetilde{h}(p, p_0))\ d\widetilde{h}(p)
\end{equation}
for any Lebesgue measurable function $\eta:(r_1,r_2)\rightarrow
[0,\infty]$ such that $\int\limits_{r_1}^{r_2}\eta(r)\,dr=1.$ }
\end{proposition}

We now state and prove the following statement.

\medskip
\begin{theorem}\label{th3}
{\sl The conclusion of Theorem~\ref{th2} holds, if instead of
conditions~(\ref{eq45}) we require that $Q\in FMO(p_0).$ }
\end{theorem}

\begin{proof} Set
$\psi(t)\,=\,\frac {1}{\left(t\,\log{\frac1t}\right)}.$
Observe that
$I(\varepsilon,
\varepsilon_0)\,:=\,\int\limits_{\varepsilon}^{\varepsilon_0}\psi(t)\,dt\,\geqslant
\log{\frac{\log{\frac{1}
{\varepsilon}}}{\log{\frac{1}{\varepsilon_0}}}}.$ Set
$\eta(t):=\psi(t)/I(\varepsilon, \varepsilon_0).$ Then, by the
relation~(\ref{eq44}), there is a constant $C> 0$ such that
$$\int\limits_{\widetilde{A}(p_0, \varepsilon, \varepsilon_0)}
Q(p)\cdot \eta^2(\widetilde{h}(p, p_0))\
d\widetilde{h}(p)=\frac{1}{I^2(\varepsilon,
\varepsilon_0)}\int\limits_{\varepsilon<\widetilde{h}(p, p_0)<
\varepsilon_0}\frac{Q(p)\, d\widetilde{h}(p)}
{\left(\widetilde{h}(p, p_0)\log\frac{1}{\widetilde{h}(p,
p_0)}\right)^2} \leqslant$$
\begin{equation}\label{eq46}
\leqslant C\cdot \left(\log{\frac{\log{\frac{1}
{\varepsilon}}}{\log{\frac{1}{\varepsilon_0}}}}\right)^{\,-1}\rightarrow
0
\end{equation}
as $\varepsilon\rightarrow 0.$ Then by~(\ref{eq10A}), taking into
account~(\ref{eq46}), conditions~(\ref{eq45}) follow, in view of
which the required conclusion follows directly from
Theorem~\ref{th2}.~$\Box$
\end{proof}

\section{Equicontinuity of families homeomorphisms}\label {sec5}

Our immediate goal is to prove the equicontinuity of the classes of
mappings consisting of Sobolev homeomorphisms of finite distortion.
First of all, let us clarify the question on the equicontinuity of
these families at the inner points of the domain. Let us recall some
definitions. Let $(X, d)$ and $\left(X^{\,{\prime}},
d^{\,{\prime}}\right)$ be metric spaces with distances $d$ and
$d^{\,{\prime}},$ respectively. A family $\frak{F}$ of mappings
$f:X\rightarrow X^{\,\prime}$ is called {\it equicontinuous at the
point} $x_0\in X,$ if for any $\varepsilon> 0$ there is $\delta> 0$
such that $d^{\,\prime}(f (x), f (x_0)) <\varepsilon$ for all $x\in
X$ such that $d(x, x_0)<\delta $ and for all $f\in\frak {F}.$ A
family $\frak {F}$ is {\it equicontinuous} if $\frak {F}$ is
equicontinuous at every point $x_0\in X.$ Everywhere below, unless
otherwise stated, $(X, d)=({\Bbb S}, \widetilde{h})$ and
$\left(X^{\,\prime}, d^{\,\prime}\right)=({\Bbb S}_*,
\widetilde{h}_*)$ are Riemannian surfaces with metrics
$\widetilde{h}$ and $\widetilde{h}_*,$ respectively.

\medskip
Let $Q:{\Bbb S}\rightarrow [0, \infty]$ be a function measurable
with respect to the measure $\widetilde{h}$ function, $Q(x)\equiv 0$
for $x\not\in D\subset {\Bbb S}.$ We say that $f:D\rightarrow {\Bbb
S}_*$ is a {\it ring $Q$-mapping at $p_0\in \overline{D},$ } if the
relation
\begin{equation}\label{eq3*!!}
M(f(\Gamma(E_1,\,E_2,\,D)))\leqslant \int\limits_A Q(p)\cdot
\eta^2(\widetilde{h}(p, p_0))\,d\widetilde{h}(p)\,,
\end{equation}
holds for some $r_0=r(p_0)>0,$ any ring
$\widetilde{A}=\widetilde{A}(p_0, r_1, r_2)=\{p\in {\Bbb S}:
r_1<\widetilde{h}(p, p_0)<r_2\},$ $0<r_1<r_2<r_0,$ and any continua
$E_1\subset \overline{\widetilde{B}(p_0, r_1)}\cap D,$ $E_2\subset
\left({\Bbb S}\setminus \widetilde{B}(p_0, r_2)\right)\cap D,$ where
$\eta:(r_1,r_2)\rightarrow [0,\infty ]$ may be arbitrary nonnegative
Lebesgue measurable function such that
\begin{equation}\label{eq28*}
\int\limits_{r_1}^{r_2}\eta(r)\,dr \geqslant 1\,.
\end{equation}

\medskip
The next definition can be found, e.g., in \cite{NP}. A domain
$D\subset {\Bbb S} $ is called a {\it uniform} if for every $r> 0$
there is $\delta> 0 $ such that $M(\Gamma(F, F ^ {\, *}, D))
\geqslant \delta $ for any continua $F$ and $F^{\, *}$ in $D,$
satisfying the conditions $\widetilde{h}(F)\geqslant r $ and
$\widetilde{h}(F^{\,*})\geqslant r.$ Domains $D_i,$ $i\in I,$ are
called {\it equi-uniform} if for each $r>0$ the above the inequality
holds for every $D_i$ with the same number $\delta.$

\medskip
For a given $\delta> 0,$ $ D\subset {\Bbb S} $ and a measurable
function $Q:D\rightarrow[0, \infty]$  with respect to the measure
$\widetilde{h},$ we denote by $\frak{R}_{Q, \delta}(D)$ the family
of all homeomorphisms $f:D\rightarrow {\Bbb S}_*\setminus G_f$ of
class $W_{\rm loc}^{1,1}$ with finite a distortion, such that
$K_f(p)\leqslant Q(p)$ for almost all $p\in D, $ where $G_f$ is some
continuum in ${\Bbb S}_*$ and
$\widetilde{h}_*(G_f)=\sup\limits_{x,y\in G_f}\widetilde{h}_*(x,
y)\geqslant \delta.$ The following statement holds.

\begin{lemma}\label{lem3}
{\sl\, The family~$\frak{R}_{Q, \delta}(D)$ is equicontinuous in
$D,$ if ${\Bbb S}_*$ is a uniform domain, $Q\in L_{\rm loc}^1$ and,
for any $p_0\in D,$ one of the following conditions is satisfied:
either~(\ref{eq45}), or $Q\in FMO(p_0).$}
\end{lemma}

\begin{proof} Since ${\Bbb S}$ is a manifold, ${\Bbb S}$ is locally
compact and locally path connected. Besides that, $f(D)$ is a domain
by Brower's theorem, see~\cite[Theorem~VI 9 and Corollary]{HW}.

Observe that the condition
\begin{equation}\label{eq1B}
\int\limits_{\varepsilon<\widetilde{h}(p,
p_0)<\varepsilon_0}Q(x)\cdot\psi^2(\widetilde{h}(p,
p_0))\,d\widetilde{h}(p)=o(I^2(\varepsilon, \varepsilon_0))
\end{equation}
holds as $\varepsilon\rightarrow 0$ for some nonnegative measurable
function $\psi:(0, \infty)\rightarrow (0, \infty)$ such that
$I(\varepsilon,
\varepsilon_0):=\int\limits_{\varepsilon}^{\varepsilon_0}\psi(t)
dt<\infty$ for some $\varepsilon_0>0$  and any $\varepsilon\in (0,
\varepsilon_0).$

Indeed, if the relations~(\ref{eq45}) hold, then~(\ref{eq1B}) holds
by Proposition~\ref{pr1A} by selecting a function $I(\varepsilon,
\varepsilon_0):=\int\limits_{\varepsilon}^{\varepsilon_0}\
\frac{dr}{\Vert Q\Vert(r)}.$ If $Q\in FMO(p_0)$ at any $p_0\in D,$
then, by the reasoning given in the proof of Theorem~\ref{th3},
conditions~(\ref{eq45}) are satisfied, and therefore, by what was
proved above, (\ref{eq1B}) holds.

Finally, by~\cite[Lemma~3.1]{RV}, the mapping $f\in \frak{R}_{Q,
\delta}(D)$ satisfies~(\ref{eq3*!!}) in $D,$ therefore the desired
conclusion follows from~\cite[Lemma~5.1]{Sev$_1$}.~$\Box$
\end{proof}

\medskip
Let us turn to the question on the equicontinuity of mappings in the
closure of a domain. For this purpose, consider the following class
mappings. Given $\delta> 0,$ $D\subset{\Bbb S},$ a continuum
$A\subset D$ and a measurable function $Q:D\rightarrow[0, \infty]$
we denote $\frak{F}_{Q, \delta, A}(D)$ the family of all Sobolev
homeomorphisms $f:D\rightarrow {\Bbb S}_*\setminus G_f$ with a
finite distortion such that $G_f\subset {\Bbb S}_*$  is some
continuum satisfying the condition
$\widetilde{h}_*(G_f)=\sup\limits_{x,y\in G_f}\widetilde{h}_*(x,
y)\geqslant \delta,$ moreover, $\widetilde{h}_*(f(A))\geqslant
\delta.$ An analogue of the following theorem was obtained
in~\cite[Theorem~3.1]{NP} for quasiconformal of mappings of the
Euclidean space.

\medskip
\begin{lemma}\label{lem4}
{\sl Let $D$ be a domain in ${\Bbb S}$ and $Q:{\Bbb S}\rightarrow
(0, \infty)$ is a function locally integrable in $D,$ $Q(x)\equiv 0$
for $x\in{\Bbb S}\setminus D.$ Assume that, for any point $p_0\in
\overline{D}$ there are $\varepsilon_0=\varepsilon_0(p_0)>0$ and a
function $\psi:(0, \infty)\rightarrow (0, \infty)$ such that
\begin{equation}\label{eq47}
I(\varepsilon,
\varepsilon_0):=\int\limits_{\varepsilon}^{\varepsilon_0}\psi(t)\,
dt<\infty\quad\forall\,\,\varepsilon\in (0, \varepsilon_0)
\end{equation}
and, in addition,
\begin{equation}\label{eq48}
\int\limits_{\varepsilon<\widetilde{h}(p,
p_0)<\varepsilon_0}Q(x)\cdot\psi^2(\widetilde{h}(p,
p_0))\,d\widetilde{h}(p)=o(I^2(\varepsilon, \varepsilon_0))\,,\quad
\varepsilon\rightarrow 0\,.
\end{equation}
Let $D_f=f(D).$ Assume also that $D$ is locally connected on
$\partial D,$ $\overline{D_f}$ is a compact in ${\Bbb S}_*$ for any
$f\in\frak{F}_{Q, \delta, A}(D),$ besides that, domains $D_f$ and
${\Bbb S}_*$ are equi-uniform over $f\in\frak{F}_{Q, \delta, A}(D).$
Then any $f\in\frak{F}_{Q, \delta, A}(D)$ has a continuous extension
$\overline{f}:\overline{D}\rightarrow \overline{D_f}$ and, besides
that, the family $\frak{F}_{Q, \delta, A}(\overline{D})$ consisting
of all extended mappings $\overline{f}:\overline{D}\rightarrow
\overline{D_f}$ is equicontinuous in $\overline{D}.$}
\end{lemma}

\begin{proof}
Observe that $\partial D_f=\partial f(D)$ is strongly accessible for
any $f\in \frak{F}_{Q, \delta, A}(D).$ Indeed, let $x_0\in
\partial D_f$ and let $U$ be an arbitrary neighborhood of $x_0.$
Choose $\varepsilon_1>0$ such that $V:=\widetilde{B}(x_0,
\varepsilon_1),$ $\overline{V}\subset U.$ Let $\partial
U\ne\varnothing$ and $\partial V\ne\varnothing.$ Now
$\varepsilon_2:=\widetilde{h}_*(\partial U,
\partial V)>0.$ Observe that, the inequalities $\widetilde{h}_*(F)\geqslant
\varepsilon_2$ and $\widetilde{h}_*(G)\geqslant \varepsilon_2$ hold
for any $F$ and $G$ in $D_f$ satisfying the conditions
$F\cap\partial U\ne\varnothing\ne F\cap\partial V$ and
$G\cap\partial U\ne\varnothing\ne G\cap\partial V.$ Now, by the
uniformity of $D_f$ there is $\delta>0$ depending only on
$\varepsilon_2$ such that
$$M(\Gamma(F, G, D_f))\geqslant \delta\,.$$
Thus, $\partial D_f$ is strongly accessible. Now,
by~\cite[Lemma~6.1]{RV} any $f\in \frak{F}_{Q, \delta, A}(D)$ has a
continuous extension $f:\overline{D}\rightarrow \overline{D_f}.$

\medskip
Observe that $\frak{F}_{Q, \delta, A}(D)\subset \frak{R}_{Q,
\delta}(D).$ Besides that, by~(\ref{eq47})--(\ref{eq48}) we obtain
that the divergence conditions~(\ref{eq45}) hold. Indeed,
in~(\ref{eq10A}) we set $\eta(t)=\psi(t)/I(\varepsilon,
\varepsilon_0)$ and let us use Proposition~\ref{pr1A}. Then the
desired conclusion immediately follows from~(\ref{eq10A}). In this
case, the equicontinuity of the family $\frak{F}_{Q, \delta, A}(D)$
in the inner points of the domain $D$ follows directly from
Lemma~\ref{lem3}.

\medskip
It remains to prove the equicontinuity of the family $\frak{F}_{Q,
\delta, A}(\overline{D})$ on $\partial D.$ Suppose the opposite.
Then there is $p_0\in
\partial D$ and a number $a>0$ such that, for each $m=1,2,\ldots$
there is a point $p_m \in \overline {D}$ and an element
$\overline{f}_m$ of the family $\overline{f}_m$$\frak{F}_{Q, \delta,
A}(\overline{D})$ such that $\widetilde{h}(p_0, p_m)<1/m$ и
$\widetilde{h}_*(\overline{f}_m(p_m), \overline{f}_m(p_0))\geqslant
a.$ Since $f_m:=\overline{f}_m|_{D}$ has a continuous extension to
the point $p_0,$ we may assume that $p_m\in D.$ In view of the same
considerations, there is a sequence
$f_m:=\overline{f}_m|_{D}$$p^{\,\prime}_m\in D,$ $p^{\,\prime}_m
\rightarrow p_0$ as $m\rightarrow\infty$ such that
$\widetilde{h}_*(f_m(p^{\,\prime}_m), \overline{f}_m(p_0))\leqslant
1/m.$ Thus
\begin{equation}\label{eq6***}
\widetilde{h}_*(f_m(p_m), f_m(p^{\,\prime}_m))\geqslant a/2\qquad
\forall\,\,m\in {\Bbb N}\,.
\end{equation}
Since $D$ is locally connected at the point $p_0\in \partial
D\subset {\Bbb S},$ and ${\Bbb S}$ is a smooth manifold, $D$ is also
locally path-connected at $p_0$ (see~\cite[Proposition~13.1]{MRSY}).
In other words, for any neighborhood $U$ of the point $p_0$ there is
a neighborhood $V\subset U$ of the same point such that $V\cap D$ is
a path-connected set. Then there is a sequence neighborhoods $V_m$
of the point $p_0$ with $\widetilde{h}(V_m)\rightarrow 0$ as
$m\rightarrow\infty,$ such that the sets $D\cap V_m $ are domains
and $D\cap V_m \subset \widetilde{B}(p_0, 2^{\,-m}).$ Without loss
of generality, passing to a subsequence, if necessary, we may assume
that $p_0\in \partial D\subset {\Bbb S},$$p_m, p^{\,\prime}_m \in
D\cap V_m.$ Join the points $p_m$ and $p^{\,\prime}_m$ of the path
$\gamma_m:[0,1]\rightarrow {\Bbb S}$ such that $\gamma_m(0)=p_m,$
$\gamma_m(1)=p^{\,\prime}_m$ and $\gamma_m(t)\in V_m\cap D$ for
$t\in (0,1).$ Denote by $C_m$ the image of the path $\gamma_m(t)$
under the mapping $f_m.$ It follows from the relation~(\ref{eq6***})
that
\begin{equation}\label{eq5.1}
\widetilde{h}_*(C_m)\geqslant a/2\qquad\forall\, m\in {\Bbb N}\,,
\end{equation}
where $\widetilde{h}_*(C_m)$ the diameter of the set $C_m$ in the
metrics $\widetilde{h}_*$ (see Figure~\ref{fig6}).
\begin{figure}[h]
\centerline{\includegraphics[scale=0.5]{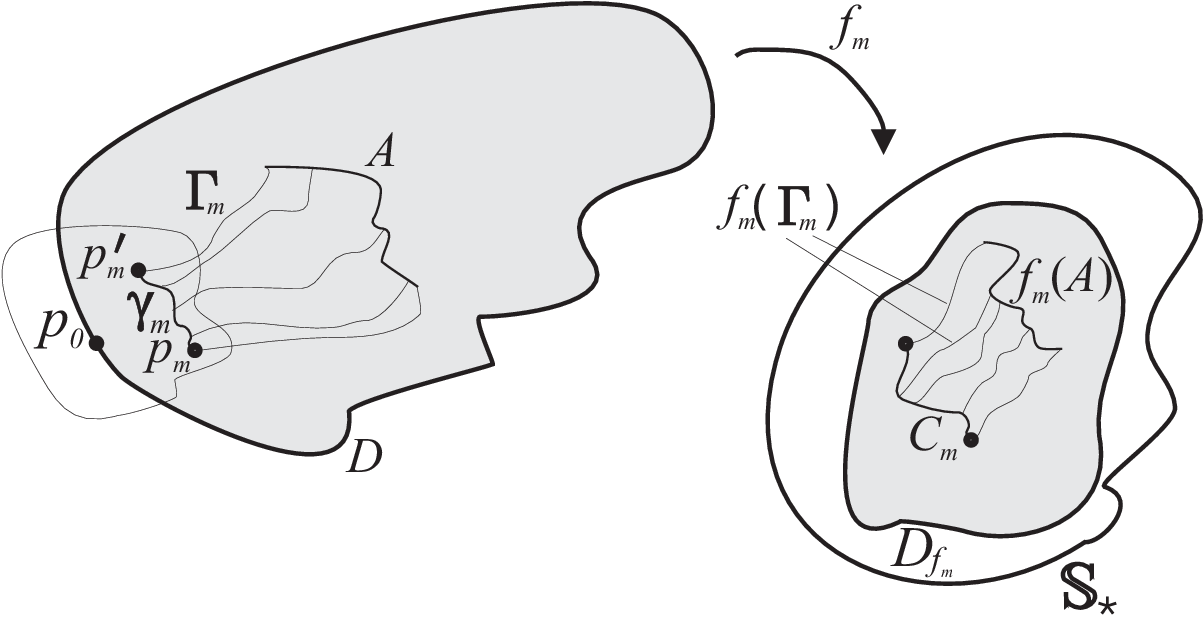}} \caption{To
the proof of Lemma~\ref{lem4}}\label{fig6}
\end{figure}
Without loss of generality, we may assume that the continuum $A$
from the definition of the family $\frak{F}_{Q, \delta, A}(D)$ is
such that $\widetilde{B}(p_0, \varepsilon_0)\cap A=\varnothing$ and
$B(p_0, 2^{\,-m})\cap A=\varnothing,$ $m=1,2,\ldots, .$ Let
$\Gamma_m$ be family of paths joining $|\gamma_m|$ and $A$ in $D.$
By~\cite[Lemma~3.1]{RV} the mapping $f\in \frak{R}_{Q, \delta}(D)$
satisfies the relation~(\ref{eq3*!!}) in $\overline{D},$ so that
\begin{equation}\label{eq10B}
M(f_m(\Gamma_m))\leqslant \int\limits_{\widetilde{A}(p_0, 2^{\,-m},
\varepsilon_0)} Q(p)\cdot \eta^2(\widetilde{h}(p,
p_0))\,d\widetilde{h}(p)
\end{equation}
for any Lebesgue measurable function $\eta:(2^{\,-m},
\varepsilon_0)\rightarrow [0, \infty]$ such that
$\int\limits_{2^{\,-m}}^{\varepsilon_0}\eta(r)\,dr \geqslant 1.$
Observe that, the function
$$\eta(t)=\left\{
\begin{array}{rr}
\psi(t)/I(2^{\,-m}, \varepsilon_0), &   t\in (2^{\,-m},
\varepsilon_0),\\
0,  &  t\in {\Bbb R}\setminus (2^{\,-m}, \varepsilon_0)\,,
\end{array}
\right. $$
где $I(\varepsilon,
\varepsilon_0):=\int\limits_{\varepsilon}^{\varepsilon_0}\psi(t)dt,$
satisfies the condition~(\ref{eq28*}) for $r_1:=2^{\,-m},$
$r_2:=\varepsilon_0,$ therefore  by~(\ref{eq47})--(\ref{eq48})
and~(\ref{eq10B}) we obtain that
\begin{equation}\label{eq11A}
M(f_m(\Gamma_m))\leqslant \alpha(2^{\,-m})\rightarrow 0
\end{equation}
as $m\rightarrow \infty,$ where $\alpha(\varepsilon)$ is some
nonnegative function converging to zero as $\varepsilon\rightarrow
0,$ which exists by~(\ref{eq47})--(\ref{eq48}).

\medskip
On the other hand, observe that $f_m(\Gamma_m)=\Gamma(C_m, f_m(A),
D_{f_m}),$ moreover, $\widetilde{h}(f_m(A))\geqslant \delta$ for any
$m\in {\Bbb N}$ by the definition of the class~$\frak{F}_{Q, \delta,
A}(D).$ Taking into account the relation~(\ref{eq5.1}) and the
definition of an equally uniform family of domains, we conclude that
there exists $\sigma> 0$ such that
$$M(f_m(\Gamma_m))=M(\Gamma(C_m, f_m(A),
D_{f_m}))\geqslant \sigma\qquad\forall\,\, m\in {\Bbb N}\,,$$
which contradicts the condition~(\ref{eq11A}). The resulting
contradiction proves the lemma.~$\Box$
\end{proof}

\medskip
By Lemma~\ref{lem4} and Proposition~\ref{pr1A} and also the
reasoning used in the proof of Theorem~\ref{th3}, we obtain the
following statement.

\begin{theorem}\label{th4}
{\sl\,Let $D$ be a domain in ${\Bbb S}$ and let $Q:{\Bbb
S}\rightarrow (0, \infty)$ be a function locally integrable in $D,$
$Q(x)\equiv 0$ for $x\in{\Bbb S}\setminus D.$ Assume that, for any
$p_0\in D,$ one of the following conditions is satisfied:
either~(\ref{eq45}), or $Q\in FMO(p_0).$
Let also $D$ be locally path-connected on $\partial D,$
$\overline{D_f}=\overline{f(D)}$ be a compactum in ${\Bbb S}_*$ for
any $f\in \frak{F}_{Q, \delta, A}(D).$ Assume that the domains $D_f$
and ${\Bbb S}_*$ are equi-uniform over $f\in\frak{F}_{Q, \delta,
A}(D).$

Then any $f\in\frak{F}_{Q, \delta, A}(D)$ has a continuous extension
$\overline{f}:\overline {D}\rightarrow \overline{D_f}$ and, in
addition, the family $\frak{F}_{Q, \delta, A}(\overline {D}),$
consisting of all extended maps $\overline{f}:
\overline{D}\rightarrow \overline{D_f},$ is equicontinuous in
$\overline {D}.$}
\end{theorem}

\section{Equicontinuity of Sobolev Classes with Branching }\label{sec6}

The question of the local and global behavior of mappings with a
branching looks much more complicated, since for mappings of
Riemannian surface estimates of the form~(\ref{eq3*!!}) have not
been established. Instead, we may only use the
estimates~(\ref{eq1A}) or~(\ref{eq15}), which are obtained in this
paper. As in the previous section, let us start by investigation of
mappings at interior points. Let us prove, first of all, the
following auxiliary statement.

\begin{proposition}\label{pr8}
{\sl\,Let $(X, d)$ be an arbitrary metric space with metric  $d$ and
$F_j,$ $j= 1,2,\ldots, $ be a sequence of continua in $X$ such that
$d(F_j)=\sup\limits_{x, y\in F_j}d(x, y)\geqslant
\delta\quad\forall\,\, j=1,2,\ldots .$
Let $x_0\in X$ and
$B(x_0, \delta/4)=\{x\in X: d(x, x_0)<\delta/4\}.$
Then there is $\varepsilon_0> 0$ and a sequence of continua $C_j$
such that $C_j \subset F_j \setminus B(x_0, \delta/4)$ and
$d(C_j)=\sup\limits_{x, y\in C_j}d(x, y)\geqslant \delta/4,\quad
\,\,j=1,2,\ldots \,.$ }
\end{proposition}

\begin{proof}
Fix $j \in {\Bbb N}.$ If $F_j\cap B(x_0, \delta/4)=\varnothing,$
there is nothing to prove. Let $F_j\cap B(x_0,
\delta/4)\ne\varnothing.$

\medskip
Since $F_j$ is a continuum in $X,$ there are $x_j, y_j \in F_j$ such
that $d(F_j)=d(x_j, y_j).$ Since $d(F_j)\geqslant\delta,$ at least
one of the points $x_j$ or $y_j$ does not belong to $B(x_0,
\delta/4), $ because, otherwise, by the triangle inequality $d(x_j,
y_j)\leqslant d(x_j, x_0)+d(x_0, y_j)<\delta/2.$ Let, for
definiteness, $x_j\in D\setminus B(x_0, \delta/4).$ Then two
situations are possible:

\medskip
1) $y_j\in B(x_0, \delta/4).$ Let $C_j$ be $x_j$-component of
$F_j\setminus B(x_0, \delta/4).$ Since $F_j$ is connected and
$F_j\cap B(x_0, \delta/4)\ne\varnothing,$ we obtain that $C_j\cap
{\overline{F_j\setminus C_j}}\ne \varnothing$ (see, e.g.,
~\cite[Section~I.5.46]{Ku}). Observe that
\begin{equation}\label{eq2F}
F_j\setminus C_j=(F_j\cap B(x_0,
\delta/4))\cup\bigcup\limits_{\alpha\in A} K_{\alpha}\,,
\end{equation}
where $A$ is some set of indices $\alpha$ and
$\bigcup\limits_{\alpha\in A} K_{\alpha}$ is a union of all
components of $F_j \setminus B(x_0, \delta/4),$ except $C_j.$
By~\cite[Theorem~1.III.46.5]{Ku}, $K_{\alpha}$ and $C_j$ are closed
disjoint sets in $F_j\setminus B (x_0, \delta/4),$ $\alpha\in A.$
Then, by~(\ref{eq2F}) the relation $C_j\cap {\overline{F_j\setminus
C_j}}\ne \varnothing$ is possible if and only if $C_j\cap
\overline{B(x_0, \delta/4)}\ne\varnothing.$ Then there is $z_j\in
C_j\cap S(x_0, \delta/4). $ By the triangle inequality
$$\delta\leqslant d(x_j, y_j)\leqslant d(x_j, z_j)+d(z_j, y_j)< d(C_j)+\delta/2\,,$$
whence it follows that $d(C_j)>\delta/2, $ as required. Consider the
second situation:

\medskip
2) $y_j\in D \setminus B(x_0, \delta/4).$ Let, as before, $C_j$ be
$x_j$-component of $F_j\setminus B(x_0, \delta/4),$ and let $D_j$
be-$y_j$-компонента $F_j\setminus B(x_0, \delta/4).$ Reasoning
similar to the above, we conclude that there are $z_j\in C_j\cap
S(x_0, \delta/4)$ and $z^{\,\prime}_j\in D_j\cap S(x_0, \delta/4).$
Then, by the triangle inequality
$$\delta\leqslant d(x_j, y_j)\leqslant d(x_j, z_j)+d(z_j, z^{\,\prime}_j)+
d(z^{\,\prime}_j, y_j)\leqslant d(C_j) + d(D_j)+ \delta/2\,,$$
whence it follows that either $d(C_j)\geqslant \delta/4,$ or
$d(D_j)\geqslant \delta/4.$ The proposition is proved.~$\Box$
\end{proof}

\medskip
The next statement concerns the situation in which the images of two
points under mappings are separated by a fixed nonzero number. It
will be shown below that in this case the length of the images of
circles centered at one of the points under these mappings is
separated from zero from below.

\begin{lemma}\label{lem1}
{\sl Let $D_{\,*}$ be a uniform domain in ${\Bbb S}_*$ such that
$\overline{D_{\,*}}$ is a compactum. Let $f_k:D\rightarrow
D_*\setminus G_k,$ $k=1,2,\ldots$ be a family of mappings open in
$D$ such that $\widetilde{h}_*(G_k)=\sup\limits_{x,y\in
G_k}\widetilde{h}_*(x, y)\geqslant \delta,$ where $G_k\subset D_*$
is some continuum and the number $\delta$ does not depend on $k.$

Suppose that $p_0\in D,$ $p_k\in D,$ $k=1,2, \ldots,$ and
$\delta_0>0$ such that $p_k\rightarrow p_0$ as $k\rightarrow\infty$
and
\begin{equation}\label{eq30B}
\widetilde{h}_*(f_k(p_k), f_k(p_0))\geqslant\delta_0\,\quad
k=1,2,\ldots\,.
\end{equation}
Then there are $l_0>0,$ $r_0>0$ and $k_0\geqslant 1$ such that
\begin{equation}\label{eq42A}
l(f_k(\widetilde{S}(p_0, r))\geqslant l_0,\quad\forall\,\,r\in
(\widetilde{h}(p_0, p_k), r_0)\,,\quad \forall\,\,k\geqslant k_0\,,
\end{equation}
where $l$ denotes the length of the path on the Riemannian surface
${\Bbb S}_*.$ }
\end{lemma}

\begin{proof}
Suppose the opposite. Then for each $i\in {\Bbb N}$ there are
$k_i>i$ and $\widetilde{h}(p_0, p_{k_i})<r_i<1/i$ such that
\begin{equation}\label{eq49}
l(f_{k_i}(\widetilde{S}(p_0, r_i)))< 1/i,\qquad i=1,2,\ldots,\qquad
r_i \rightarrow 0
\end{equation}
as $i\rightarrow\infty.$
Without loss of generality, we may assume that the sequence numbers
$k_i,$ $i=1,2, \ldots $ is increasing. Let $\zeta_i,$ $i=1,2,
\ldots, $ be an arbitrary sequence of points from
$f_{k_i}(\widetilde{S}(p_0, r_i)).$ Since $\overline{D_{\,*}}$ is a
compactum in ${\Bbb S}_*,$ we may assume that $\zeta_i\rightarrow
\zeta_0$ as $i\rightarrow\infty,$ $\zeta_0\in \overline{D_{\,*}}.$
Note that $\zeta_i=f_{k_i}(p^{\,\prime}_i),$ $p^{\,\prime}_i\in
\widetilde{S}(p_0, r_i),$ and that
\begin{equation}\label{eq2E}
G_{k_i}\subset D_*\setminus \overline{f_{k_i}(\widetilde{B}(p_0,
r_i))}\,,
\end{equation}
because by the condition $G_k\subset D_*\setminus f_k(D)$ for any
$k\in {\Bbb N}.$ Since $f_{k_i}$  is open for any $i\in {\Bbb N},$
we obtain that
\begin{equation}\label{eq7B}
\partial f_{k_i}(\widetilde{B}(p_0,
r_i)) \subset f_{k_i}(\widetilde{S}(p_0, r_i))\,.
\end{equation}
The assumption~(\ref{eq49}) implies that
$$\widetilde{h}_*(f_{k_i}(\widetilde{S}(p_0,
r_i)))\rightarrow 0$$
as $i\rightarrow\infty,$ $\widetilde{h}_*(f_{k_i}(\widetilde{S}(p_0,
r_i))):=\sup\limits_{p_*, q_*\in f_{k_i}(\widetilde{S}(p_0,
r_i))}\widetilde{h}_*(p_*, q_*).$
Now, for any $s\in {\Bbb N}$ there exists a number $i_s\in {\Bbb N}$
such that
\begin{equation}\label{eq50}
f_{k_i}(\widetilde{S}(p_0, r_i))\subset \widetilde{B}(\zeta_0,
1/s)\,,\qquad i\geqslant i_s\,.
\end{equation}
By Proposition~\ref{pr8}, there is $s_0\in {\Bbb N} $ and a
sequence~$E_{k_i}$ of continua such that
\begin{equation}\label{eq51}
E_{k_i}\subset G_{k_i}\setminus \widetilde{B}(\zeta_0, 1/s_0),\quad
\widetilde{h}_*(E_{k_i})\geqslant \delta/4\,,\quad i=1,2,\ldots \,.
\end{equation}
We fix $s>s_0$ and consider the family
$\Gamma(f_{k_i}(\widetilde{S}(p_0, r_i)), E_{k_i}, D_*)$ for
$i\geqslant i_s,$ see Figure~\ref{fig5}.
\begin{figure}[h]
\centerline{\includegraphics[scale=0.5]{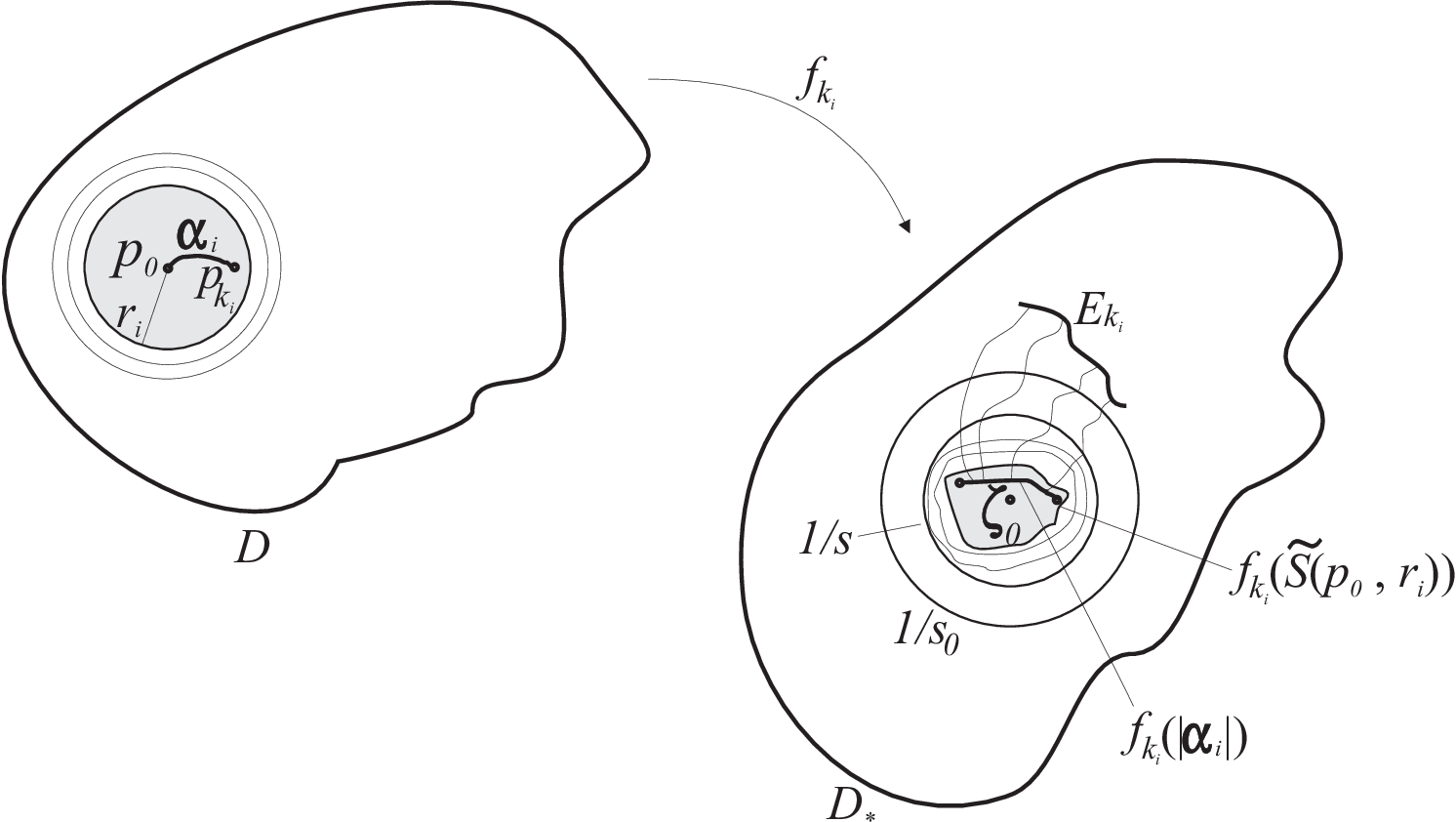}} \caption{To
the proof of Lemma~\ref{lem1}}\label{fig5}
\end{figure}
Let $\gamma\in \Gamma(f_{k_i}(\widetilde{S}(p_0, r_i)), E_{k_i},
D_*),$ i.e.,  $\gamma=\gamma(t),$ $t\in (0, 1),$ $\gamma(0)\in
f_{k_i}(\widetilde{S}(p_0, r_i)),$ $\gamma(1)\in E_{k_i}$ and
$\gamma(t)\in D_*$ for $t\in (0, 1).$ By~(\ref{eq50})
and~(\ref{eq51}) $|\gamma|\cap \widetilde{B}(\zeta_0,
1/s)\ne\varnothing\ne |\gamma|\cap
(D_*\setminus\widetilde{B}(\zeta_0, 1/s)),$ therefore
by~\cite[Theorem~1.I.5, $\S\,46$]{Ku} there is $t_1\in (0, 1)$ such
that $\gamma(t_1)\in\widetilde{S}(\zeta_0, 1/s).$ Without loss of
generality we may assume that $\gamma(t)\in D_*\setminus
\overline{\widetilde{B}(\zeta_0, 1/s)}$ for $t>t_1.$ Set
$\gamma_1:=\gamma|_{[t_1, 1]}.$ Again, by ~(\ref{eq50})
and~(\ref{eq51}) $|\gamma_1|\cap \widetilde{B}(\zeta_0,
1/s_0)\ne\varnothing\ne |\gamma_1|\cap
(D_*\setminus\widetilde{B}(\zeta_0, 1/s_0)),$ therefore
by~\cite[Theorem~1.I.5, $\S\,46$]{Ku} there is $t_2\in (t_1, 1):$
$\gamma_1(t_2)\in\widetilde{S}(\zeta_0, 1/s_0).$ Without loss of
generality we may assume that $\gamma_1(t)\in \widetilde{B}(\zeta_0,
1/s_0)$ for $t\in (t_1, t_2).$ Set $\gamma_2:=\gamma_1|_{[t_1,
t_2]},$ $\gamma_2\in \Gamma(\widetilde{S}(\zeta_0, 1/s),
\widetilde{S}(\zeta_0, 1/s_0), \widetilde{A}(\zeta_0, 1/s, 1/s_0)),$
$\widetilde{A}(\zeta_0, 1/s, 1/s_0))=\{p_*\in {\Bbb S}_*:
1/s<\widetilde{h}_*(p_*, \zeta_0)<1/s_0\}.$ From the above it
follows that
$$\Gamma(f_{k_i}(\widetilde{S}(p_0, r_i)), E_{k_i}, D_*)>\Gamma(\widetilde{S}(\zeta_0, 1/s),
\widetilde{S}(\zeta_0, 1/s_0), \widetilde{A}(\zeta_0, 1/s,
1/s_0))\,,\quad i\geqslant i_s\,,$$
and, therefore, in view of~\cite[Theorem~1(c)]{Fu} and by
Proposition~\ref{pr5}
$$M(\Gamma(f_{k_i}(\widetilde{S}(p_0, r_i)), E_{k_i}, D_*))\leqslant$$
\begin{equation}\label{eq52}
\leqslant M(\Gamma(\widetilde{S}(\zeta_0, 1/s),
\widetilde{S}(\zeta_0, 1/s_0), \widetilde{A}(\zeta_0, 1/s,
1/s_0)))\rightarrow 0
\end{equation}
as $s\rightarrow \infty$ for $i \geqslant i_s.$ Let us fix
$\varepsilon> 0$ and find for it a number $S=S(\varepsilon)$ such
that
$$M(\Gamma(\widetilde{S}(\zeta_0, 1/s), \widetilde{S}(\zeta_0,
1/s_0), \widetilde{A}(\zeta_0, 1/s, 1/s_0)))<\varepsilon, \quad
s>S(\varepsilon)\,.$$ Set
$I_0=I_0(\varepsilon):=i_{S(\varepsilon)}.$ Now, it follows
from~(\ref{eq52}) that
\begin{equation}\label{eq53}
M(\Gamma(f_{k_i}(\widetilde{S}(p_0, r_i)), E_{k_i},
D_*))<\varepsilon, \quad i>I_0=I_0(\varepsilon)\,.
\end{equation}
Since ${\Bbb S}$ is a smooth manifold, we may consider that the
balls $\widetilde{B}(p_0, r_i)$ are path connected for any $i\in
{\Bbb N}.$ Let $\alpha_i$ be a path joining $p_{k_i}$ and $p_0$ in
$\widetilde{B}(p_0, r_i).$ By~(\ref{eq30B})
$\widetilde{h}_*(f_{k_i}(|\alpha_i|))\geqslant \delta_0.$ Now, by
the definition of a uniform domain
\begin{equation}\label{eq54}
M(\Gamma(f_{k_i}(|\alpha_i|), E_{k_i}, D_*))>\varepsilon_1\,,\quad
\forall\,\,i\in {\Bbb N}\,.
\end{equation}
On the other hand, by~(\ref{eq2E}) and (\ref{eq7B})
$$\Gamma(f_{k_i}(|\alpha_i|), E_{k_i},
D_*)>\Gamma(f_{k_i}(\widetilde{S}(p_0, r_i)), E_{k_i}, D_*)\,,$$
whence by~\cite[Theorem~1(c)]{Fu} and also by~(\ref{eq54})
\begin{equation}\label{eq55}
\varepsilon_1<M(\Gamma(f_{k_i}(|\alpha_i|), E_{k_i}, D_*))\leqslant
M(\Gamma(f_{k_i}(\widetilde{S}(p_0, r_i)), E_{k_i}, D_*))\,.
\end{equation}
The inequalities~(\ref{eq55}) and (\ref {eq53}) contradict each
other, which proves~\ref{eq42A}).~$\Box$
\end{proof}

\medskip
Let $D\subset {\Bbb S}$ and $D_*\subset {\Bbb S}_*$ be fixed
domains. Given $\delta>0$ and a measurable function $Q:{\Bbb
S}\rightarrow (0, \infty)$ with a respect to the measure
$\widetilde{h},$ $Q(x)\equiv 0$ for $x\in{\Bbb S}\setminus D,$ we
demote by $\frak{G}_{Q, \delta}(D)$ the family of all open mappings
$f:D\rightarrow D_*\setminus G_f$ satisfying the
relation~(\ref{eq1A}) in $D$ such that
$\widetilde{h}_*(G_f)=\sup\limits_{x,y\in G_f}\widetilde{h}_*(x,
y)\geqslant \delta,$ where $G_f\subset D_*$ is some continuum.

\begin{lemma}\label{lem5}
{\sl Assume that $Q$ satisfies~(\ref{eq45}) in $D,$ or $Q\in
FMO(p_0)$ at any point $p_0\in D.$ If $D_*$ is uniform and
$\overline{D_{\,*}}$ is a compactum in ${\Bbb S}_*,$ then the family
$\frak{G}_{Q, \delta}(D)$ is equicontinuous in $D.$}
\end{lemma}

\begin{proof}
In view of the reasoning used in the proof of Theorem~\ref{th3}, it
suffices to establish Lemma~\ref{lem5} in the case when $Q$
satisfies relations~(\ref{eq45}) in $D.$

Suppose that the conclusion of the lemma does not hold. Then there
are $p_0\in D,$ $p_k\in D,$ $k=1,2, \ldots,$ $f_k\in\frak{G}_{Q,
\delta}(D)$ and $\delta_0$ such that $p_k \rightarrow p_0$ as
$k\rightarrow\infty$ and
\begin{equation}\label{eq30A}
\widetilde{h}_*(f_k(p_k), f_k(p_0))\geqslant\delta_0\,.
\end{equation}
By Lemma~\ref{lem1}, there are $l_0>0$ and $r_0>0$ such that
\begin{equation}\label{eq42B}
l(f_k(\widetilde{S}(p_0, r))\geqslant l_0
\end{equation}
for all $r\in (\widetilde{h}(p_0, p_k), r_0),$ for any $k\geqslant
k_0$ and some $k_0\geqslant 1.$ Without loss of generality, we may
assume that $r_0<\varepsilon_0, $ where $\varepsilon_0 $ is a number
from~(\ref{eq45}) that exists by the condition of the lemma.
In this case, the function $$\rho(p)= \left \{\begin{array}{rr}
1/l_0\ , & \ p\in D_{\,*}\ ,
\\ 0\ ,  &  p\not\in D_{\,*}
\end{array} \right.$$
is admissible for $\Gamma_k^{r_0},$ consisting of the union of the
paths $f_k(\widetilde{S}(p_0, r))$ over all $r\in(\widetilde{h}(p_0,
p_k), r_0),$ $k=1,2, \ldots .$
In this case, by the definition of the modulus of families of paths,
\begin{equation}\label{eq57}
M(\Gamma_k^{r_0})\leqslant (1/l_0^2)\cdot
\widetilde{h_*}(D_{\,*})<\infty\,,
\end{equation}
because $\overline{D_{\,*}}$ is a compactum in ${\Bbb S}_*.$ On the
other hand, by Lemma~\ref{lem4A} and also by~(\ref{eq45}) we obtain
that
\begin{equation}\label{eq56}
M(\Gamma_k^{r_0})\geqslant \int\limits_{\widetilde{h}(p_0,
p_k)}^{r_0} \frac{dr}{\Vert\,Q\Vert(r)}\rightarrow\infty
\end{equation}
as $k\rightarrow\infty.$
The relations~(\ref{eq57}) and (\ref{eq56}) contradict each other,
which refutes the assumption made in~(\ref{eq30A}).~$\Box$
\end{proof}

\medskip
As usual, we formulate the main results of this section for Sobolev
classes. For this purpose, consider the following definition of the
family of mappings. Given numbers $\delta> 0,$ $N\in {\Bbb N},$ a
domain $D\subset {\Bbb S}$ and a function $Q:{\Bbb S}\rightarrow (0,
\infty)$ that is measurable with respect to the measure
$\widetilde{h},$ $Q(x)\equiv 0$ for $x\in{\Bbb S}\setminus D,$
denote by $\frak{S}_{Q, \delta, N}(D)$ a family of all open discrete
mappings $f:D\rightarrow D_*\setminus G_f$ of the class $W_{\rm
loc}^{1, 1}(D)$ with a finite distortion such that $K_f(p)\leqslant
Q(p)$ for almost all $p\in D,$ $N(f, D)\leqslant N$ and
$\widetilde{h}_*(G_f)=\sup\limits_{x,y\in G_f}\widetilde{h}_*(x,
y)\geqslant \delta,$ where $G_f\subset D_*$ is some continuum. The
following theorem holds.

\begin{theorem}\label{th5}
{\sl Suppose that the function $Q$ satisfies the
relations~(\ref{eq45}) in $D,$ or $Q\in FMO(p_0)$ at any point
$p_0\in D.$ If the domain $D_*$ is uniform and $\overline{D_{\,*}}$
is a compactum in ${\Bbb S}_*,$ then the family $\frak{S}_{Q,
\delta, N}(D)$ is equicontinuous at any point $p_0\in D.$}
\end{theorem}

\medskip
{\it Proof} immediately follows from Lemma~\ref{lem5} and
Theorem~\ref{thOS4.2}.~$\Box$

\medskip
Let us turn to the study of equicontinuity at the points of the
boundary. First of all, similarly to Lemma~\ref{lem1}, we prove the
following statement.

\begin{lemma}\label{lem8}
{\sl Let $D$ and $D _ *$ be domains in ${\Bbb S}$ and ${\Bbb S}_ *,$
respectively, $p_0\in \partial D,$ $p_k\in D,$ $k=1,2, \ldots,$
$p_k\rightarrow p_0$ as $k\rightarrow\infty,$ and let a domain $D$
has a locally quasiconformal boundary. Let
$f_k:\overline{D}\rightarrow \overline{D_*},$ $k=1,2,\ldots$ be a
family of mappings such that $f_k|_D$ is open and closed in $D.$
Suppose that

\medskip
1) the domains $D_k:=f_k (D)$ and $D_*$ are equi-uniformly over
$k\in{\Bbb N},$ in addition, $\overline{D_{\, *}}$ is a compactum in
${\Bbb S}_*; $

\medskip
2) there is a number $\delta>0$ with the following property: for any
$k\in{\Bbb N}$ there is a continuum $A_k\subset f_k(D),$ such that
$\widetilde{h}_*(A_k)\geqslant \delta>0,$ moreover,
$\widetilde{h}(f_k^{\,-1}(A_k), \partial D)\geqslant \delta>0;$

\medskip
3) there is $\delta_0>0$ such that
\begin{equation}\label{eq30C}
\widetilde{h}_*(f_k(p_k), f_k(p_0))\geqslant\delta_0\,\quad
\forall\,\,k=1,2,\ldots\,.
\end{equation}
Then there are $l_0> 0,$ $r_0>0$ and $k_0\geqslant 1$ such that
\begin{equation}\label{eq42C}
l(f_k(\widetilde{S}(p_0, r)\cap D)\geqslant l_0,\quad\forall\,\,r\in
(\widetilde{h}(p_0, p_k), r_0)\,,\quad \forall\,\,k\geqslant k_0\,,
\end{equation}
where $l$ denotes the length of the dashed line $\widetilde{S}(p_0,
r)\cap D$ on the Riemannian surface ${\Bbb S}_*.$  }
\end{lemma}

\begin{proof}
Suppose the opposite. Then for each $i\in{\Bbb N} $ there are
$k_i>i$ and $\widetilde{h}(p_0, p_ {k_i})<r_i<1 / i,$ for which
\begin{equation}\label{eq49A}
l(g_i(\widetilde{S}(p_0, r_i)\cap D))< 1/i,\qquad
i=1,2,\ldots,\qquad r_i\rightarrow 0\,,
\end{equation}
as $i\rightarrow\infty,$ $g_i:=f_{k_i}.$ Without loss of generality,
we may assume that the sequence numbers $k_i,$ $i=1,2,\ldots$ is
increasing. We also denote $p_{k_i}:=q_i.$

Since $D$ is a domain with a locally quasiconformal boundary, by
Lemma~\ref{pr7} there is a sequence of neighborhoods $U_m$ of the
point $p_0,$ $m=1,2,\ldots ,$ contracting to $p_0$ such that
$U_m\cap D$ is connected and, moreover, for each $m=1,2,\ldots ,$
there is a cut $\gamma_m$ of $D$ such that $U_m\cap D \subset
D\setminus|\gamma_m|,$ $|\gamma_m|\subset \widetilde{S}(p_0,
r_{i_m})$ for some subsequence $r_{i_m}$ of the sequence $r_i,$ $i=
1,2, \ldots. $ Since $f_k$ is continuous in $\overline {D},$ there
is a sequence $p_k^{\,\prime}\in D$ such that
$\widetilde{h}_*(f_k(p_0), f_k(p_k^{\,\prime}))<1/k.$ Put
$q_i^{\,\prime}:=p_{k_i}^{\,\prime}.$ Then it follows
from~(\ref{eq30C}) that
\begin{equation}\label{eq30D}
\widetilde{h}_*(g_i(q_i), g_i(q_i^{\,\prime}))\geqslant\delta_0/2
\end{equation}
for $i\geqslant i_0\in {\Bbb N}.$ For any $m\in {\Bbb N},$ we may
find $j_m\in {\Bbb N},$ $j_m>i_0,$ such that $q_{j_m},
q_{j_m}^{\,\prime}\in U_m.$ Now, by~(\ref{eq30D}) we obtain that
\begin{equation}\label{eq30E}
\widetilde{h}_*(g_{j_m}(q_{j_m}),
g_{j_m}(q_{j_m}^{\,\prime}))\geqslant\delta_0/2\,\quad
\forall\,\,m=1,2,\ldots\,.
\end{equation}

\medskip
Let us show that
\begin{equation}\label{eq2G}
A_{j_m}\subset D_*\setminus g_{j_m}(U_m\cap \overline{D})\qquad
\forall\,\,m\geqslant m_0\,.
\end{equation}
Indeed, if $y_l\in A_{j_{m_l}}\cap g_{j_{m_l}}(U_{m_l}\cap
\overline{D})$ for arbitrarily large $l=1,2,\ldots$ and some
increasing sequence $m_l,$ $l=1,2,\ldots,$ then
$y_l=g_{j_{m_l}}(x_l),$ $x_l\in U_{m_l}\cap \overline{D}$ and, at
the same time, $x_l\in g^{\,-1}_{j_{m_l}}(A_{j_{m_l}}).$ Since
$U_{m_l}$ contracts to $p_0\in \partial D, $ then $x_l\rightarrow
p_0$ as $l\rightarrow\infty,$ which contradicts the condition
$\widetilde{h}(f_k^{\,-1}(A_k),
\partial D)\geqslant \delta>0,$ $g_{j_{m_l}}=f_{k_{j_{m_l}}}.$ Thus, (\ref{eq2G}) holds.

\medskip
Now let $d_m$ be the component of $D\setminus|\gamma_m|,$ containing
$U_m\cap D.$ Let us show that
\begin{equation}\label{eq7C}
\partial g_{j_m}(d_m)\cap g_{j_m}(D)\subset
g_{j_m}(|\gamma_m|)
\end{equation}
for any $m=1,2,\ldots .$

We fix $m\in{\Bbb N}$ and consider $y_m\in \partial g_{j_m}(d_m)\cap
g_{j_m}(D).$ Then there is $y_{mk}\in g_{j_m}(d_m),$
$y_{mk}\rightarrow y_m$ as $k\rightarrow\infty.$ Since $f_k$ are
open, $g_{j_m}(D)$ is a domain, so we may assume that $y_{mk}\in
g_{j_m}(d_m)\cap g_{j_m}(D).$ Since  $y_{mk}\in g_{j_m}(d_m),$ then
$y_{mk}=g_{j_m}(\eta_{mk}),$ $\eta_{mk}\in d_m.$ By lemma~\ref{pr7}
$d_n\subset U_m $ for all $n\geqslant n(M),$ therefore $d_m$ also
contract to the point $p_0$ as $m\rightarrow \infty. $ Thus, we may
assume that $\overline{d_m} $ is a compactum in ${\Bbb S}, $ and
that $\eta_{mk}\rightarrow \eta_0$ as $k\rightarrow\infty.$ Observe
that the case $\eta_0 \in \partial D $ is impossible, because now
$y_m\in C(g_{j_m},
\partial D)\subset
\partial g_{j_m}(D)$ by the closeness of the mapping $g_{j_m}=f_{k_{j_m}},$
which contradicts the choice of $y_m.$ Then $\eta_0 \in D.$ Two
situations are possible: 1) $\eta_0\in d_m$ and 2)
$\eta_0\in|\gamma_m|.$ Observe that the case 1) is impossible,
because now $g_{j_m}(\eta_0)=y_m$ and $y_m$ is an inner point of the
set $g_{j_m}(d_m)$ by the openness of the mapping $g_{j_m},$ which
also contradicts the choice of $y_m.$ Thus, the
inclusion~(\ref{eq7C}) is proved.

\medskip
The further course of reasoning largely repeats the scheme of the
proof of Lemma~\ref{lem1}. Let $\xi_m,$ $m=1,2,\ldots, $ be an
arbitrary sequence of points from $|\gamma_m|.$ Since
$\overline{D_{\,*}}$ is a compactum in ${\Bbb S}_*,$ without loss of
generality, we may assume that $\zeta_m:=g_{j_m}(\xi_m)\rightarrow
\zeta_0$ as $m\rightarrow\infty,$ $\zeta_0\in \overline{D_{\,*}}.$
It follows from~(\ref{eq49A}) that
$$\widetilde{h}_*(g_{j_m}(|\gamma_m|))\rightarrow 0$$
as $m\rightarrow\infty,$
$\widetilde{h}_*(g_{j_m}(|\gamma_m|)):=\sup\limits_{p_*, q_*\in
g_{j_m}(|\gamma_m|)}\widetilde{h}_*(p_*, q_*).$
Now, for any $s\in {\Bbb N}$ there is a number $m_s\in {\Bbb N}$
such that
\begin{equation}\label{eq50A}
g_{j_m}(|\gamma_m|)\subset \widetilde{B}(\zeta_0, 1/s)\,,\qquad
m\geqslant m_s\,.
\end{equation}
By Proposition~\ref{pr8} there is $s_0\in {\Bbb N}$ and a sequence
of continua~$E_{j_m}$ such that
\begin{equation}\label{eq51A}
E_{j_m}\subset A_{j_m}\setminus \widetilde{B}(\zeta_0, 1/s_0),\quad
\widetilde{h}_*(E_{j_m})\geqslant \delta/4\,,\quad i=1,2,\ldots \,.
\end{equation}
We fix $s>s_0$ and consider the family $\Gamma(g_{j_m}(|\gamma_m|),
E_{j_m}, g_{j_m}(D))$ for $m\geqslant m_s.$ Let $\gamma\in
\Gamma(g_{j_m}(|\gamma_m|), E_{j_m}, g_{j_m}(D)),$ i.e.,
$\gamma=\gamma(t),$ $t\in (0, 1),$ $\gamma(0)\in
f_{k_i}(g_{j_m}(|\gamma_m|),$ $\gamma(1)\in E_{j_m}$ and
$\gamma(t)\in g_{j_m}(D)$ for $t\in (0, 1).$ By~(\ref{eq50A})
and~(\ref{eq51A}) $|\gamma|\cap \widetilde{B}(\zeta_0,
1/s)\ne\varnothing\ne |\gamma|\cap
(g_{j_m}(D)\setminus\widetilde{B}(\zeta_0, 1/s)),$ therefore,
according to~\cite[Theorem~1.I.5,$\S\,46$]{Ku} there is $t_1\in (0,
1)$ such that $\gamma(t_1)\in\widetilde{S}(\zeta_0, 1/s).$ Without
loss of generality, we may assume that $\gamma(t)\in
g_{j_m}(D)\setminus \overline{\widetilde{B}(\zeta_0, 1/s)}$ for
$t>t_1.$ Put $\gamma_1:=\gamma|_{[t_1, 1]}.$ Again, by~(\ref{eq50A})
and~(\ref{eq51A})  $|\gamma_1|\cap \widetilde{B}(\zeta_0,
1/s_0)\ne\varnothing\ne |\gamma_1|\cap
(g_{j_m}(D)\setminus\widetilde{B}(\zeta_0, 1/s_0)),$ therefore
by~\cite[Theorem~1.I.5,$\S\,46$]{Ku} there is $t_2\in(t_1, 1)$ such
that $\gamma_1(t_2)\in \widetilde{S}(\zeta_0, 1/s_0).$ Without loss
of generality, we may assume that $\gamma_1(t)\in
\widetilde{B}(\zeta_0, 1/s_0)$ for $t\in (t_1, t_2).$ Put
$\gamma_2:=\gamma_1|_{[t_1, t_2]},$ $\gamma_2\in
\Gamma(\widetilde{S}(\zeta_0, 1/s), \widetilde{S}(\zeta_0, 1/s_0),
\widetilde{A}(\zeta_0, 1/s, 1/s_0)),$ $\widetilde{A}(\zeta_0, 1/s,
1/s_0))=\{p_*\in {\Bbb S}_*: 1/s<\widetilde{h}_*(p_*,
\zeta_0)<1/s_0\},$ see Figure~\ref{fig7}.
\begin{figure}[h]
\centerline{\includegraphics[scale=0.5]{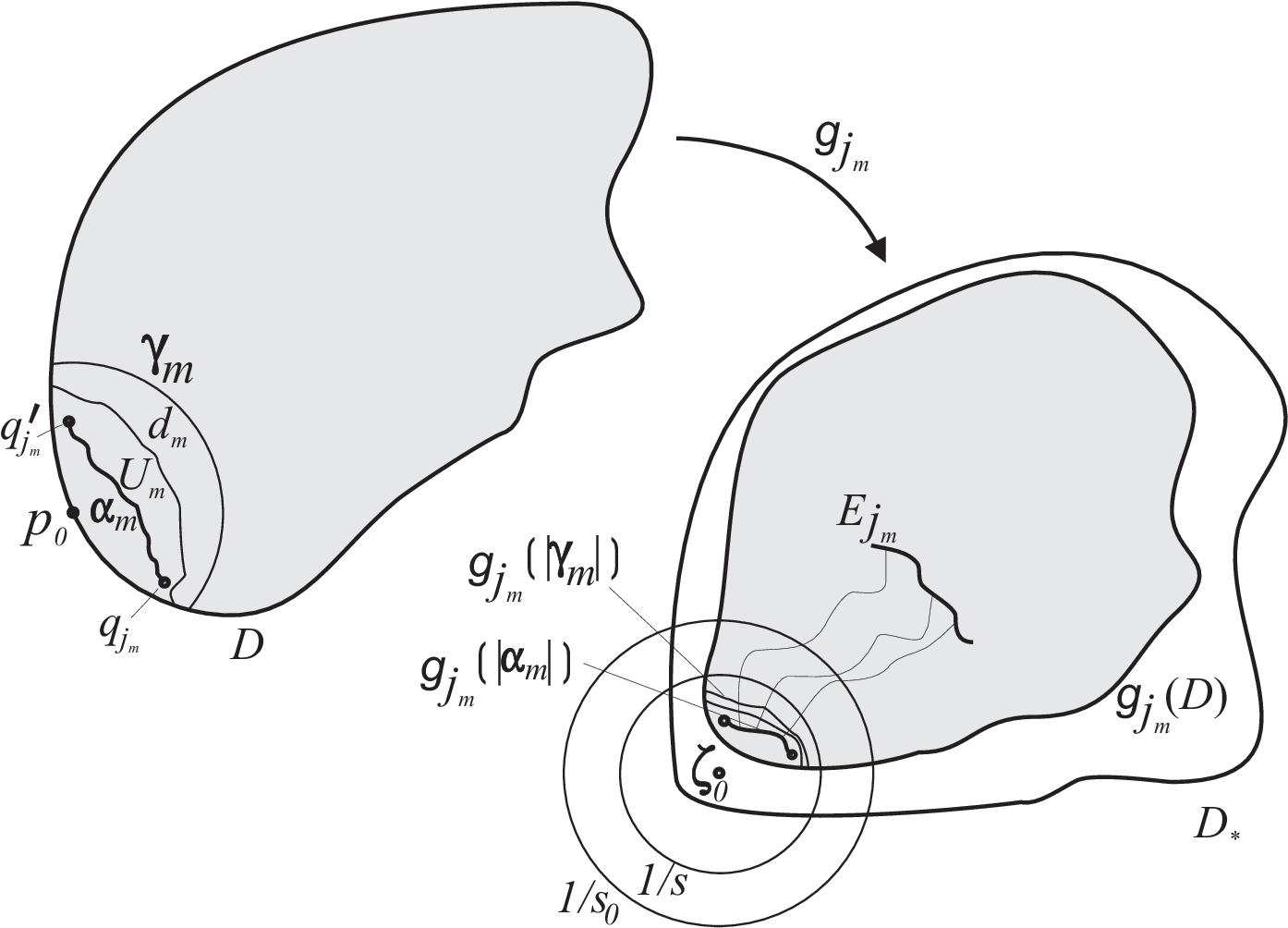}} \caption{To
the proof of Lemma~\ref{lem8}}\label{fig7}
\end{figure}
Hence it follows that
$$\Gamma(g_{j_m}(|\gamma_m|), E_{j_m},
g_{j_m}(D))>\Gamma(\widetilde{S}(\zeta_0, 1/s),
\widetilde{S}(\zeta_0, 1/s_0), \widetilde{A}(\zeta_0, 1/s, 1/s_0))$$
for $m\geqslant m_s.$ Thus, by~\cite[Theorem~1(c)]{Fu} and by
Proposition~\ref{pr5}
$$M(\Gamma(g_{j_m}(|\gamma_m|), E_{j_m}, g_{j_m}(D)))\leqslant$$
\begin{equation}\label{eq52A}
\leqslant M(\Gamma(\widetilde{S}(\zeta_0, 1/s),
\widetilde{S}(\zeta_0, 1/s_0), \widetilde{A}(\zeta_0, 1/s,
1/s_0)))\,,
\end{equation}
$$M(\Gamma(\widetilde{S}(\zeta_0, 1/s), \widetilde{S}(\zeta_0, 1/s_0),
\widetilde{A}(\zeta_0, 1/s, 1/s_0)))\rightarrow 0$$
as $s\rightarrow \infty.$
Given $\varepsilon>0$  we may find a number $S=S(\varepsilon)$ such
that
$$M(\Gamma(\widetilde{S}(\zeta_0, 1/s), \widetilde{S}(\zeta_0,
1/s_0), \widetilde{A}(\zeta_0, 1/s, 1/s_0)))<\varepsilon, \quad
s>S(\varepsilon)\,.$$ Set
$M_0=M_0(\varepsilon):=m_{S(\varepsilon)}.$ It follows
from~(\ref{eq52A}) that
\begin{equation}\label{eq53A}
M(\Gamma(g_{j_m}(|\gamma_m|), E_{j_m}, g_{j_m}(D)))<\varepsilon,
\quad m>M_0=M_0(\varepsilon)\,.
\end{equation}
Now let us join the points $q_{j_m}$ and $q_{j_m}^{\,\prime}$ by the
path $\alpha_m\subset U_m.$ By assumption~(\ref{eq30D})
$\widetilde{h}_*(g_{j_m}(|\alpha_m|))\geqslant \delta_0/2.$ Then by
the definition of a uniform domain, there exists $\varepsilon_1> 0$
such that
\begin{equation}\label{eq54A}
M(\Gamma(g_{j_m}(|\alpha_m|), E_{j_m},
g_{j_m}(D)))>\varepsilon_1\,,\qquad m\in {\Bbb N}\,.
\end{equation}
On the other hand, by~(\ref{eq2G}) and (\ref{eq7C})
$$\Gamma(g_{j_m}(|\alpha_m|), E_{j_m},
g_{j_m}(D))>\Gamma(g_{j_m}(|\gamma_m|), E_{j_m}, g_{j_m}(D))\,,$$
whence by~\cite[Theorem1(c)]{Fu} and also by~(\ref{eq54A}) we obtain
that
\begin{equation}\label{eq55A}
\varepsilon_1<M(\Gamma(g_{j_m}(|\alpha_m|), E_{j_m},
g_{j_m}(D)))\leqslant M(\Gamma(g_{j_m}(|\gamma_m|), E_{j_m},
g_{j_m}(D)))
\end{equation}
for $m=1,2,\ldots .$ The inequalities~(\ref{eq55A}) and
(\ref{eq53A}) contradict each other, which
proves~(\ref{eq42C}).~$\Box$
\end{proof}

\medskip
Let$D$ and $D_*$ be domains in ${\Bbb S}$ and ${\Bbb S}_*,$
respectively. Given $\delta>0$ and a measurable function $Q:{\Bbb
S}\rightarrow (0, \infty)$ with respect the measure $\widetilde{h},$
$Q(x)\equiv 0$ for $x\in{\Bbb S}\setminus D,$ denote by
$\frak{L}_{Q, \delta}(D)$ the family of all open discrete mappings
$f:D\rightarrow D_*$ with~(\ref{eq1A}) for which: 1) there is a
continuum $G_f\subset D_*$ such that $f:D\rightarrow D_*\setminus
G_f$ и $\widetilde{h}_*(G_f):=\sup\limits_{x,y\in
G_f}\widetilde{h}_*(x, y)\geqslant \delta;$ 2) there exists a
continuum $A_f\subset f(D)$ such that
$\widetilde{h}_*(A_f)\geqslant\delta$ and
$\widetilde{h}(f^{\,-1}(A_f),
\partial D)\geqslant\delta.$ The following theorem holds.

\begin{theorem}\label{th6}
{\sl\,Assume that the following conditions are satisfied:

\medskip
1) $Q$ is locally integrable in $D$ and either $Q$
satisfies~(\ref{eq45}), or $Q\in FMO(p_0)$ at any point $p_0\in
\overline{D};$

\medskip
2) domains $D_f=f(D)$ and $D_*$ are equi-uniform over $f\in
\frak{L}_{Q, \delta}(D);$

\medskip
3) a domain $D$ has a locally quasiconformal boundary;

\medskip
4) the set $\overline{D_{\,*}}$ is a compactum in ${\Bbb S}_*.$

\medskip
Then any mapping $f\in\frak{L}_{Q, \delta}(D)$ has a continuous
extension $\overline{f}:\overline{D}\rightarrow \overline{D_{\,*}}$
and a family $\frak{L}_{Q, \delta}(\overline{D})$ of all extended
mappings $\overline{f}$ is equicontinuous in $\overline{D}.$}
\end{theorem}

\begin{proof}
Arguing similarly to the beginning of the proof of Lemma~\ref{lem4},
we conclude that the domain $D_f=f(D),$ $f\in \frak{L}_{Q,
\delta}(D),$ has a strongly accessible boundary. In this case, the
possibility of continuous extension
$\overline{f}:\overline{D}\rightarrow \overline{D_{\,*}}$ follows
from Theorem~\ref{th1}, and the equicontinuity of $\frak{L}_{Q,
\delta}(\overline{D})$ in $D$ is from Lemma~\ref{lem5}, since
$\frak{L}_{Q, \delta}(D)\subset \frak{G}_{Q, \delta}(D).$ It remains
to show the equicontinuity of the family $\frak{L}_{Q,
\delta}(\overline{D})$ in $\partial D.$

\medskip
Suppose the opposite. Then there are $p_0\in \partial D,$ $p_k\in
\overline{D},$ $k=1,2, \ldots,$ $f_k\in\frak{L}_{Q,
\delta}(\overline{D})$ and $\delta_0$ such that $p_k\rightarrow p_0$
as $k\rightarrow\infty$ and
\begin{equation}\label{eq30F}
\widetilde{h}_*(f_k(p_k), f_k(p_0))\geqslant\delta_0\,.
\end{equation}
By Theorem~\ref{th1} we may assume that $p_k\in D,$ besides that, by
Lemma~\ref{lem1} there are $l_0>0$ and $r_0>0$ such that for some
$k_0\geqslant 1$
\begin{equation}\label{eq42D}
l(f_k(\widetilde{S}(p_0, r))\geqslant l_0,\quad\forall\,\,r\in
(\widetilde{h}(p_0, p_k), r_0)\,,\quad \forall\,\,k\geqslant k_0\,,
\end{equation}
where $l$ denotes the length of a path on the Riemannian surface
${\Bbb S}_*.$ Without loss of generality, we may assume that
$r_0<\varepsilon_0,$ where $\varepsilon_0$ is a number
from~(\ref{eq45}), existing by the conditions of the lemma.
In this case, the function $$\rho(p)= \left \{\begin{array}{rr}
1/l_0\ , & \ p\in D_{\,*}\ ,
\\ 0\ ,  &  p\not\in D_{\,*}
\end{array} \right.$$
is admissible for $\Gamma_k^{r_0},$ consisting from the union of all
dished lines $f_k(\widetilde{S}(p_0, r)),$ $k=1,2, \ldots ,$ over
$r\in(\widetilde{h}(p_0, p_k), r_0).$
In this case, by the definition of the modulus of families of paths,
we obtain that
\begin{equation}\label{eq57A}
M(\Gamma_k^{r_0})\leqslant (1/l_0^2)\cdot
\widetilde{h_*}(D_{\,*})<\infty\,,
\end{equation}
because $\overline{D_{\,*}}$ is a compactum in ${\Bbb S}_*.$ On the
other hand, by Lemma~\ref{lem4A} and by the conditions~(\ref{eq45})
we have that
\begin{equation}\label{eq56A}
M(\Gamma_k^{r_0})\geqslant \int\limits_{\widetilde{h}(p_0,
p_k)}^{r_0} \frac{dr}{\Vert\,Q\Vert(r)}\rightarrow\infty
\end{equation}
as $k\rightarrow\infty.$ The relations~(\ref{eq57A}) and
(\ref{eq56A}) contradict each other,, which refutes the assumption
made in~(\ref{eq30F}).~$\Box$
\end{proof}

\medskip
Let $D$ and $D_*$ be domains in ${\Bbb S}$ and ${\Bbb S}_*,$
respectively. Given $\delta>0,$ a natural number $N\geqslant 1$ and
a function $Q:{\Bbb S}\rightarrow (0, \infty),$ $Q(x)\equiv 0$ for
$x\in{\Bbb S}\setminus D,$ measurable with respect to the measure
$\widetilde{h},$ denote by $\frak{M}_{Q, \delta, N}(D)$ the family
of all open discrete and closed mappings $f:D\rightarrow D_*$ of the
class $W_{\rm loc}^{1, 1}(D)$ with a finite distortion for which: 1)
there is a continuum $G_f\subset D_*$ such that $f:D\rightarrow
D_*\setminus G_f$ and $\widetilde{h}_*(G_f):=\sup\limits_{x,y\in
G_f}\widetilde{h}_*(x, y)\geqslant \delta;$ 2) there is a continuum
$A_f\subset f(D)$ such that $\widetilde{h}_*(A_f)\geqslant\delta$
and $\widetilde{h}(f^{\,-1}(A_f),
\partial D)\geqslant\delta;$ 3) $K_f(p)\leqslant Q(p)$ for any $p\in D;$
4) $N(f, D)\leqslant N.$

\medskip
The following theorem holds.

\begin{theorem}\label{th7}
{\sl\,Assume that the following conditions are satisfied:

\medskip
1) $Q$ is locally integrable in $D$ and either $Q$
satisfies~(\ref{eq45}), or $Q\in FMO(p_0)$ at any point $p_0\in
\overline{D};$

\medskip
2) domains $D_f=f(D)$ and $D_*$ are equi-uniform over $f\in
\frak{M}_{Q, \delta, N}(D);$

\medskip
3) a domain $D$ has a locally quasiconformal boundary;

\medskip
4) the set $\overline{D_{\,*}}$ is a compactum in ${\Bbb S}_*.$

\medskip
Then any mapping $f\in\frak{M}_{Q, \delta, N}(D)$ has a continuous
extension $\overline{f}:\overline{D}\rightarrow \overline{D_{\,*}}$
and a family $\frak{M}_{Q, \delta, N}(\overline{D})$ of all extended
mappings $\overline{f}$ is equicontinuous in $\overline{D}.$ }
\end{theorem}

\medskip
{\it Proof if Theorem~\ref{th7}} immediately follows from
Theorems~\ref{th6} and \ref{thOS4.2}.~$\Box$

%=================Список литературы====================
%\end{fulltext}

\medskip
\medskip
{\bf \noindent Evgeny Sevost'yanov} \\
{\bf 1.} Zhytomyr Ivan Franko State University,  \\
40 Bol'shaya Berdichevskaya Str., 10 008  Zhytomyr, UKRAINE \\
{\bf 2.} Institute of Applied Mathematics and Mechanics\\
of NAS of Ukraine, \\
1 Dobrovol'skogo Str., 84 100 Slov'yans'k,  UKRAINE\\
esevostyanov2009@gmail.com

\end{document}